\documentclass{article}
\usepackage{amsfonts}
\usepackage{amsmath,amsthm}

\usepackage[all]{xy}
\usepackage{graphicx}

\usepackage{amssymb}
\usepackage{latexsym}

\DeclareMathAlphabet\EuFrak{U}{euf}{m}{n}	
\SetMathAlphabet\EuFrak{bold}{U}{euf}{b}{n}	

\newcommand{\plihat}{{\wa \psi}_l^i}
\newcommand{\pmjhat}{{\wa \psi}_m^j}
\newcommand{\aot}{a \otimes_\mu t}

\newcommand{\ra}{\rightarrow}
\newcommand{\lra}{\longrightarrow}

\newcommand{\hra}{\hookrightarrow}

\newcommand{\ovl}{\overline}

\newcommand{\wa}{\widehat}
\newcommand{\wt}{\widetilde}

\newcommand{\sC}{{\it C*}-}
\newcommand{\bC}{{\mathbb C}}

\newcommand{\bT}{{\mathbb T}}

\newcommand{\bZ}{{\mathbb Z}}
\newcommand{\bM}{{\mathbb M}}
\newcommand{\bN}{{\mathbb N}}
\newcommand{\bP}{{\mathbb P}}
\newcommand{\bS}{{\mathbb S}}

\newcommand{\ud}{{{\mathbb U}(d)}}
\newcommand{\sud}{{{\mathbb {SU}}(d)}}

\newcommand{\eps}{\varepsilon}

\newcommand{\mA}{\mathcal A}
\newcommand{\mB}{\mathcal B}

\newcommand{\mC}{\mathcal C}
\newcommand{\mE}{\mathcal E}
\newcommand{\mF}{\mathcal F}
\newcommand{\mG}{\mathcal G}

\newcommand{\mL}{\mathcal L}
\newcommand{\mM}{\mathcal M}

\newcommand{\mO}{\mathcal O}

\newcommand{\mR}{\mathcal R}

\newcommand{\mZ}{\mathcal Z}

\newcommand{\tend}{{\bf end}\mA}

\newcommand{\zro}{{C(X^\rho)}}
\newcommand{\spzro}{ X^\rho }

\newcommand{\soro}{{\mO_{\rho , \eps}}}

\newcommand{\ii}{\iota,\iota}

\newcommand{\mrs}{\mM^r,\mM^s}
\newcommand{\ers}{\mE^r,\mE^s}
\newcommand{\rhors}{\rho^r , \rho^s}

\newcommand{\sgrs}{\sigma_G^r,\sigma_G^s}

\newcommand{\wE}{\wa {\mE}}

\newcommand{\mSUE}{ {{\bf {SU}} \mE} }
\newcommand{\mUE}{ { {\bf U} \mE} }

\newcommand{\mSG}{ SG }

\newcommand{\cpen}{ \mA \rtimes_\rho^\mE \bN }

\newcommand{\m}{\mA \odot_\mu {^0\mathcal O_{\mE}}}
\newcommand{\mc}{\mA \rtimes_\mu \wa G }
\newcommand{\mo}{\mA_x \otimes_{\mu_x} {\mO_d}}

\newcommand{\coe}{\mO_{\mE}}

\newcommand{\cog}{\mO_G}
\newcommand{\cmg}{\mO^{G}_{\mM}}
\newcommand{\aoe}{{^0\mathcal O_{\mE}}}
\newcommand{\aog}{{^0\mathcal O_G}}

\newtheorem{thm}{Theorem}[section]
\newtheorem{cor}[thm]{Corollary}
\newtheorem{lem}[thm]{Lemma}
\newtheorem{prop}[thm]{Proposition}

\newtheorem{defn}[thm]{Definition}

\theoremstyle{definition}
\newtheorem{ex}{Example}[section]

\theoremstyle{remark}
\newtheorem{rem}{Remark}[section]

\numberwithin{equation}{section}

\begin{document}

\author{{\sf Ezio Vasselli}
                         \\{\it Dipartimento di Matematica}
                         \\{\it University of Rome "La Sapienza"}
			 \\{\it P.le Aldo Moro, 2 - 00185 Roma - Italy }
                         \\{\sf vasselli@mat.uniroma2.it}}

\title{ Crossed products by endomorphisms,\\vector bundles\\and\\group duality, II}
\maketitle

\begin{abstract}
We study \sC algebra endomorphims which are special in a weaker sense w.r.t. the notion introduced by Doplicher and Roberts. We assign to such endomorphisms a geometrical invariant, representing a cohomological obstruction for them to be special in the usual sense. Moreover, we construct the crossed product of a \sC algebra by the action of the dual of a (nonabelian, noncompact) group of vector bundle automorphisms. These crossed products supply a class of examples for such generalized special endomorphisms.

\bigskip

\noindent {\em AMS Subj. Class.:} 46L05, 46L08, 22D35.

\noindent {\em Keywords:} Crossed Products; Tensor \sC categories; $C_0(X)$-algebras; Vector bundles.

\end{abstract}


\section{Introduction.}
\label{intro}

Let $\mA$ be a unital \sC algebra with centre $\mZ$. For each pair of endomorphisms $\rho,\sigma$, we consider the intertwiners spaces $(\rho,\sigma) := \left\{ t \in \mA : \sigma (a) t = t \rho (a) , a \in \mA \right\}$. Every $(\rho,\sigma)$ has an obvious structure of Banach $\mZ$-bimodule $z,t \mapsto zt$, with $z \in \mZ$, $t \in (\rho,\sigma)$. In particular, if $\rho$ is the identity $\iota$, we get a Hilbert $\mZ$-bimodule, endowed with the natural $\mZ$-valued scalar product $t,t' \mapsto t^* t'$, $t , t' \in (\iota , \sigma)$. We say that an endomorphism $\rho$ is {\it inner} if $(\iota,\rho)$ is non-zero and finitely generated as a right Hilbert $\mZ$-module, i.e. if there is a finite set $\left\{ \psi_l \right\} \subset (\iota,\rho)$ such that $\rho (1) = \sum_l \psi_l \psi_l^*$. Our terminology is justified by the fact that if $\rho$ is inner then, for $a \in \mA$,
\begin{equation}
\label{12} 
\rho (a) = \sum \limits_l \rho (a) \psi_l \psi_l^* = \sum \limits_l \psi_l a \psi_l^* \ ,
\end{equation}
\noindent so that $\rho$ is induced by every set of generators of $(\iota,\rho)$, which play a role analogous to the unitaries for inner automorphisms. Viceversa, every finitely generated Hilbert $\mZ$-bimodule contained in $\mA$ in the above sense induces an inner endomorphism of $\mA$, {\em via} (\ref{12}). The previous notion of inner endomorphism has been considered and used in \cite{DR88}, in the case in which $\mA$ has trivial centre; the terminology adopted in that context is {\it Hilbert space in a C*-algebra}, since the scalar product is $\bC$-valued.

As for \sC algebra automorphisms, one can define the crossed product by an endomorphism, in such a way that the endomorphism itself becomes induced by an isometry (a one dimensional Hilbert space, in our terminology; see \cite{Pas80,Cun77,Exe00}). More general constructions have been done by using Hilbert spaces of isometries, as in \cite{DR89A}, and after in \cite{Sta93}. In the first-cited reference, the crossed product of a \sC algebra by an endomorphism is related to the duality theory for compact groups, in the following sense: the set $\tend$ of unital endomorphisms of $\mA$, with arrows the intertwiners spaces $(\rho,\sigma)$, has a natural structure of a tensor category when endowed with the tensor structure
\begin{equation}
\label{def_enda}
\left\{ 
\begin{array}{l}
\rho,\sigma \mapsto \rho \sigma := \rho \circ \sigma  \\ 
t,t' \mapsto t \times t' := t \rho (t') = \rho'(t')t \in (\rho \sigma , \rho' \sigma')
\end{array}
\right.
\end{equation}
\noindent where $t \in (\rho,\rho') , t' \in (\sigma,\sigma')$. Thus, one could ask which are the conditions such that some full subcategory of ${\bf end} \mA$ is isomorphic to a group dual. The basic step for the answer is given (in the case of a single endomorphism) by the following definitions:

\begin{defn}[Permutation symmetry]
\label{def11}
A unital endomorphism $\rho$ of $\mA$ has permutation symmetry if there is a unitary representation $p \mapsto \eps (p)$ of the group $\bP_\infty$ of finite permutations of $\bN$ in $\mA$, such that:
\begin{equation}\label{ps1} \eps (\bS p) = \rho \circ \eps (p) \end{equation}
\begin{equation}\label{ps2} \eps := \eps (1,1) \in (\rho^2 , \rho^2) \end{equation}
\begin{equation}\label{ps3} \eps (s,1) \ t = \rho (t) \ \eps (r,1) \ \ , \ \ t \in (\rhors) \ , \end{equation}
\noindent where $(r,s) \in \bP_{r+s}$ permutes the first $r$ terms with the remaining $s$, and $\bS$ is the shift $(\bS p)(1) := 1$, $(\bS p)(n) := 1 + p(n-1)$, $p \in \bP_\infty$ ($\bP_n$, $n \in \bN$, denotes the usual permutation group of $n$ objects). The above properties imply that
\begin{equation}
\label{weak_perm}
\eps(p) \in (\rho^n , \rho^n) \ \ , \ \ n \in \bN \ , \ p \in \bP_n \ ;
\end{equation}
\noindent we say that $\rho$ has {\bf weak permutation symmetry} if just (\ref{ps1}),(\ref{weak_perm}) hold.
\end{defn}

\begin{defn}[Special Conjugate Property]
\label{def12}
Let $\rho$ be an endomorphism of $\mA$ with permutation symmetry. $\rho$ has dimension $d > 1$ and satisfies the special conjugate property if there exists an isometry $R \in (\iota , \rho^d)$ such that

\begin{enumerate}
\item $R^* \rho (R) = (-1)^{d-1} d^{-1} 1 \ ,$
\item $RR^* =  P_{\eps,d} := \frac 1{d!} \sum \limits_{p \in \bP_d} \mathrm{sign}(p) \eps (p)$.
\end{enumerate}
\end{defn}

\noindent The previous definitions model at an algebraic level respectively the symmetry of the tensor product and the existence of a determinant for a representation of a compact group. An endomorphism satisfying Def.\ref{def12} is called {\em special}, according to the terminology of \cite[\S 4]{DR89A}. Let us now denote by $\wa \rho$ the full tensor subcategory of $\tend$ with objects $\rho^r$, $r \in \bN$, and arrows $(\rhors)$. In this setting, it can be proven that if  $\mA$ is a unital \sC algebra with trivial centre, and $\rho \in \tend$ a special endomorphism, then $\wa \rho$ is isomorphic to the category of tensor powers of the defining representation of a compact Lie group $G \subseteq \sud$ (\cite[Thm.4.1,Lemma 4.6]{DR89A}). The above-cited theorem is an important step towards the Doplicher-Roberts duality for compact groups, and is motivated by questions arising in the context of algebraic quantum field theory (see \cite{DR90} and related references). It has also been applied to the nontrivial centre case in \cite{BL01,BL04}, by maintaining (in essence) the properties Def.\ref{def11}, Def.\ref{def12}; the result (in the case of a single endomorphism) is a crossed product where the endomorphism becomes induced by a free module of the type $\bC^d \otimes \mZ$, and a group $G \subseteq \sud$ is recovered as in the trivial-centre case.

As pointed out before, in general the natural objects that provide the notion of inner endomorphism are (not free) Hilbert \sC bimodules over the centre of the given \sC algebra i.e., by the Serre-Swan Theorem, modules of continuous sections of vector bundles. In the present paper, which is a continuation of \cite{Vas04}, we study the analogues of the previously quoted constructions in a generalized setting: we consider a nontrivial centre  and relax permutation symmetry and special conjugate property, which do not hold in the examples that we consider (shift endomorphisms in the setting of Cuntz-Pimsner algebras). Instead of compact Lie groups, we will deal with certain non-compact groups of vector bundle automorphisms.

The present paper is organized as follows.

\

In Section \ref{dual_act} we go through the same line of \cite{DR89A}, and generalize the notion of dual action of a Lie group on a \sC algebra $\mA$. Let $\mE \ra X$ be a vector bundle, $\mUE$ the unitary group of $\mE$, $G \subseteq \mUE$ a closed group. A dual $G$-action is a functor $\mu : \wa G \ra \tend$, where $\wa G$ is the category of tensor powers of the defining representation $u_0$ of $G$ over $\mE$ (see (\ref{def_g_dual})). In this setting, we construct a crossed product \sC algebra $\mc$, containing the module $\wE$ of continuous sections of $\mE$ (Thm.\ref{thm38}). $\wE$ induces on $\mA$ the endomorphism $\mu (u_0) \in \tend$, as by (\ref{12}). There is a strongly continuous action $G \ra {\bf aut}(\mc)$, in such a way that (with some technical assumptions) the fixed-point algebra is isomorphic to $\mA$ (Lemma \ref{lem38}). The crossed product $\mc$ supplies a generalization of the notion of Hilbert \sC system considered in \cite{BL97,BL01,BL04}.

\

In Section \ref{spec_end} we study \sC algebra endomorphisms which are {\em{quasi-special}} (resp. {\em weakly special}) according to generalized (less restrictive) notions of permutation symmetry (Def.\ref{def_gps}) and special conjugate property (Def.\ref{def13}). In particular, we allow the existence of intertwiners that do not satisfy (\ref{ps3}); otherwise, we will talk about {\em symmetry intertwiners} (Def.\ref{def_strict_int}). Let $\rho \in {\bf end}\mA$ be a (unital) weakly special endomorphism; we define
\begin{equation}
\label{def_zro}
\zro := \left\{ 
f \in \mA \cap \mA' : \rho (f) = f
\right\} \ ,
\end{equation}
\noindent and assign to $\rho$ a geometrical invariant $c(\rho) := c_0 (\rho) \oplus c_1 (\rho)$, belonging to the sheaf cohomology group 
\begin{equation}
\label{def_hzro}
H^{0,2}( \spzro , \bZ ) := H^0( \spzro , \bZ ) \oplus H^2( \spzro , \bZ )
\end{equation}
\noindent (Def.\ref{cro}). The first summand in $H^0( \spzro , \bZ )$ corresponds to a generalized dimension of the permutation symmetry; the second summand in $H^2( \spzro , \bZ )$ denotes a cohomological obstruction to get the special conjugate property in the sense of Def.\ref{def12}. Special endomorphisms in the sense of \cite{DR89A,BL04} thus correspond to elements of the type $d \oplus 0 \in H^{0,2}( \spzro , \bZ )$, $d \in \bN$. We give duality results when a dual action is defined, characterizing $\wa \rho$ as the category of tensor powers of a suitable $G$-Hilbert bimodule (Prop.\ref{cor_dual_1}, Prop.\ref{rhoine}).

\

Section \ref{cp_spec} contains our main result. We prove that if $\rho$ is a weakly special endomorphism in the sense of Def.\ref{cro}, then for every vector bundle $\mE \ra \spzro$ with rank $c_0(\rho)$ and first Chern class $c_1(\rho)$, a natural dual $\mSUE$-action is induced on $\mA$ (Thm.\ref{sue_action}; here $\mSUE \subset \mUE$ denotes the group of unitaries with determinant $1$). Moreover, the \sC algebra $\soro \subseteq \mA$ generated by symmetry intertwiners is characterized as a \sC algebra bundle, and a unique family of compact groups $\left\{ G_x \subseteq \sud \right\}_{x \in \spzro}$ is associated to $\rho$, in such a way that the fibre of $\soro$ over $x \in \spzro$ is the fixed-point algebra of the Cuntz algebra $\mO_d$ w.r.t. the canonical $G_x$-action.

\section{Keywords.}
\label{key}

\subsection{Some definitions and references.}
The Doplicher-Roberts theory for {\em special endomorphisms} can be found in \cite{DR87,DR89A}.

Let $X$ be a locally compact Hausdorff space. A {\em $C_0(X)$-algebra} is a \sC algebra $\mA$ endowed with a nondegenerate inclusion $C_0(X) \hra ZM(\mA)$, where $ZM(\mA)$ is the centre of the multiplier algebra (\cite{Kas88}). In the unital case (with $X$ compact), then $C(X)$ is just a unital subalgebra of $\mA \cap \mA'$. We will denote by ${\bf aut}_X \mA$ (resp. ${\bf end}_X \mA )$ the set of automorphisms (resp. endomorphisms) of $\mA$ preserving the $C_0(X)$-algebra structure. Every $C_0(X)$-algebra can be regarded as an upper-semicontinuous bundle of \sC algebras (\cite{Bla96,Nil96}); in particular, every continuous bundle of \sC algebras over $X$ (\cite[Chp.10]{Dix},\cite{KW95}) is a $C_0(X)$-algebra.

For general properties of Cuntz-Pimsner algebras, we refere to \cite{Pim97}. In the present paper, we will follow the approach used in \cite{DPZ97}. We refere to the last-cited reference also for the notion of {\em Hilbert bimodule in a \sC algebra} considered in the previous introduction.

For basic properties of {\em vector bundles}, we refere to \cite{Ati,Kar}. About the Cuntz-Pimsner algebra of a vector bundle, and the associated canonical shift endomorphism, we refere to \cite{Vas04}; for readers which are unfamiliar with the above-cited reference, we recommend the following subsection.

\subsection{Tensor categories in the Cuntz-Pimsner algebra of a vector bundle.}
\label{apdx}

Let $X$ be a compact Hausdorff space, $\mE \ra X$ a rank $d$ vector bundle, $d \in \bN$. It is well-known that the Hilbert $C(X)$-bimodule $\wE$ of continuous sections of $\mE$ is finitely generated; we consider a finite set of generators $\left\{ \psi_l \right\}_{l=1}^n \subset \wE$, $n \in \bN$. The Cuntz-Pimsner algebra $\coe$ associated with $\wE$ can be described in terms of generators and relations:
\begin{equation}
\label{genrel}
\psi_l^* \psi_m = \left \langle \psi_l , \psi_m \right \rangle  \ \ , \ \ 
\psi_l f = f \psi_l \ \ , \ \ 
\sum_l \psi_l \psi_l^* = 1 \ \ ,
\end{equation}
\noindent where $f \in C(X)$, and $\left \langle \cdot , \cdot \right \rangle$ denotes the $C(X)$-valued scalar product of $\wE$. In the sequel, we will regard at $\wE$ as a (generating) subset of $\coe$. It turns out that $\coe$ is a \sC algebra bundle over $X$, with fibre the Cuntz algebra $\mO_d$ (\cite[Prop.2]{Vas04a},\cite{Rob}). Let $r,s \in \bN$; we denote by $\mE^r$ the $r$-fold tensor power of $\mE$, and by $(\ers)$ the Banach $C(X)$-bimodule of vector bundle morphisms from $\mE^r$ into $\mE^s$. By the Serre-Swan theorem, $(\ers)$ can be identified with the set of $C(X)$-bimodule operators from $\wE^r$ into $\wE^s$, where $\wE^r$ denotes the $r$-fold (internal) tensor power of $\wE$. In particular, we define $\iota := \mE^0 := X \times \bC$, so that there is a natural identification $\wE = (\iota , \mE)$. Let $\psi , \psi' , \varphi \in \wE$, and $\theta_{\psi' , \psi} \in (\mE , \mE)$ be the operator $\theta_{\psi' , \psi} \varphi := \psi' \left \langle \psi , \varphi \right \rangle$; every $(\ers)$ can be regarded as a Banach subspace of $\coe$, and by (\ref{genrel}) we find that the following relations hold in $\coe$:
\[
\left \langle \psi , \psi' \right \rangle = \psi^* \psi' 
\ \ , \ \
\theta_{\psi' , \psi} = \psi' \psi^* \ \ .
\]
\noindent If $L := \left\{ l_1 , \ldots , l_r \right\}$ is a multiindex with lenght $r \in \bN$, we introduce the notation
\begin{equation}
\label{def_psiL}
\psi_L := \psi_{l_1} \cdots \psi_{l_r} \in (\iota , \mE^r) \subset \coe
\end{equation}
\noindent so that 
$(\ers) = {\mathrm{span}} 
\left\{ 
f \psi_M \psi_L^* \ , \ f \in C(X) , |L| = r , |M| = s
\right\}$.
Let $\mUE$ denote the unitary group of $\mE$, and $G \subseteq \mUE$ a closed group. Then, $G$ acts in a natural way on $\coe$ by $C(X)$-automorphisms: if $g \in G$, the associated automorphism $\wa g \in {\bf aut}_X \coe$ is defined by
\begin{equation}
\label{24} 
\wa{g} (t) := g_s \cdot t \cdot g^*_r  \in (\ers) 
\ \ , \ \
t \in (\ers) \ ,
\end{equation}
\noindent where $g_r := g \otimes \cdots \otimes g \in (\mE^r , \mE^r)$ denotes the $r$-fold tensor power. We denote by $\wa G$ the category with objects $\mE^r$, $r \in \bN$, and arrows
\begin{equation}
\label{def_g_dual}
(\ers)_G := \left\{ t \in (\ers) : \wa g (t) = t \right\} \ .
\end{equation}
\noindent It turns out that $\wa G$ is a tensor \sC category in the sense of \cite{DR89,LR97}, with $(\ii)_G = (\ii) = C(X)$. We also denote by $\cog \subseteq \coe$ the \sC subalgebra generated by the set $\left\{ (\ers)_G , r,s \in \bN \right\}$. It is clear that $\cog$ is a bundle; we denote by $(\cog)_x \subseteq \mO_d$, $x \in X$, the corresponding fibres. Let us now consider the canonical endomorphism
\begin{equation}
\label{def_ce}
\sigma \in {\bf end}_X \coe \ \ , \ \ \sigma (t) := 1 \otimes t 
\ \ , \ \ 
t \in (\ers) \ \ .
\end{equation}
\noindent It turns out that $\sigma$ is inner in the sense of the previous section, with $(\iota , \sigma) = (\iota , \mE) = \wE$ (\cite[Prop.4.2]{Vas}). Moreover, $\sigma$ has permutation symmetry, induced by the flip operator 
\begin{equation}
\label{defsim}
\theta = \sum \limits_{l,m} {\psi}_m {\psi}_l {\psi}^*_m {\psi}^*_l 
\in (\mE^2 , \mE^2)_{\mUE} 
\ : \
\theta \psi \psi' = \psi' \psi \  .
\end{equation}
\noindent Thus, we get unitary representations $\bP_r \ni p \mapsto \theta (p) \in (\mE^r , \mE^r)$ (\cite[Rem.4.5]{Vas04}). $\sigma$ satisfies the special conjugate property Def.\ref{def12} if and only if the first Chern class of $\mE$ vanishes (\cite[Lemma 4.2]{Vas04}): in general, instead of the partial isometry appearing in Def.\ref{def12}, we have the Hilbert $C(X)$-bimodule $\mR := (\iota , \lambda \mE) \subset (\iota , \sigma^d)$, corresponding to the module of continuous sections of the totally antisymmetric line bundle $\lambda \mE \subset \mE^d$. If $\left\{ R_i \right\} \subset (\iota , \sigma^d)$ is a finite set of generators for $\mR$, then we find
\begin{equation}
\label{eq_P} 
\sum \limits_i R_i R_i^{*} = P_{\theta,d} := 
\frac 1{d!} \sum \limits_{p \in \bP_d} \mathrm{sign} (p) \theta (p) \ ;
\end{equation}
\noindent moreover,
\begin{equation}
\label{eq_scp1}
R^* \sigma (R') = (-1)^{d-1} d^{-1} R^* R' \ \ , \ \ R , R' \in \mR \ .
\end{equation}
\noindent For every $G \subseteq \mSUE$, it turns out $\mR \subseteq (\iota , \mE^d)_G$. Since $\cog$ is $\sigma$-stable, the endomorphism $\sigma_G := \sigma |_{\cog} \in {\bf end}_X \cog$ is well defined; it turns out that $(\ers)_G = (\sgrs)$, $r,s \in \bN$, so that there is an isomorphism of tensor \sC categories $\wa G \simeq \wa \sigma_G$ (\cite[Cor.4.4]{Vas04}). We recall a last result: let $x \in X$, and
\begin{equation}
\label{eq_def_G^x}
G^x := \left\{ 
               u \in \ud \ : \ 
               \wa u (y) = y \ , \ 
               y \in (\cog)_x \subseteq \mO_d 
       \right\} \ .
\end{equation}
\noindent If $G \subseteq \mSUE$, then $(\cog)_x = \mO_{G^x}$, $x \in X$. A suitable topology (compatible with the one of $\mUE$) can be defined on the set $\left\{ G^x \right\}_{x \in X}$, in such a way to get a bundle $\mG \ra X$. It turns out that the stabilizer ${\bf aut}_{\cog} \coe$ of $\cog$ in $\coe$ is isomorphic to the group $\mSG$ of continuous sections of $\mG$. Moreover, there are inclusions $G \subseteq \mSG \subseteq \mUE$, and $\wa G = \wa{\mSG}$, i.e. $(\ers)_G = (\ers)_{\mSG}$, $r,s \in \bN$. In general, $G$ may not coincide with $\mSG$ (\cite[\S 4.3]{Vas04}).

%
%
%
%
%
%

\section{Crossed Products by Dual Actions.}
\label{dual_act}

\subsection{Basic Definitions.}

\begin{defn}
\label{defn21}
Let $\mA$ be a unital \sC algebra, $X$ a compact Hausdorff space, $\mE \ra X$ a vector bundle with a closed group $G \subseteq \mUE$. A dual $G$-action on $\mA$ is an injective functor of tensor \sC categories $\mu : {\wa G} \ra \tend$.
\end{defn}

$\wa G$ is the tensor \sC category defined by (\ref{def_g_dual}); since $(\ers)_G = (\ers)_{\mSG}$, $r,s \in \bN$, {\em we may always assume that} $G = \mSG$ (see Sec.\ref{apdx}). In the sequel we will denote $\rho := \mu (\mE) \in \tend$, so that a dual $G$-action (all the objects of ${\wa G}$ being tensor powers of $\mE$) is actually an injective functor of tensor \sC categories $\mu : {\wa G} \ra \wa \rho$, where $\wa \rho$ is the full subcategory of $\tend$ generated by the powers of $\rho$. Thus, the condition that $\mu$ preserves the tensor product is expressed by
\[
\mu ( y \otimes y' ) = \mu (y) \times \mu (y') = \mu (y) \cdot \rho^r \circ \mu (y')
\ ,
\]
\noindent $y \in ( \ers )_G$, $y' \in (\mE^{r'} , \mE^{s'})_G$ (see (\ref{def_enda})). In particular, $\mu (f) = \mu (1 \otimes f) = \rho \circ \mu (f)$, $f \in C(X) = (\ii)_G$, where $1 \in (\mE , \mE)$ denotes the identity map; thus, $\mu (f) \in \zro$ (see (\ref{def_zro})), and we have a unital monomorphism $\mu : C(X) \hra \zro$. So that, a surjective map $\mu_* : \spzro \ra X$ is induced, and it is not restrictive to assume that $X = \spzro$: in fact, we may replace $\mE$ with the pullback $\mE_\rho := \mE \times_X \spzro \ra \spzro$, and identify $G$ as a subgroup of $\mUE_\rho$.

Thus, $\cog$ and $\mA$ have an obvious structure of $\zro$-algebra. Since $(\sgrs) = (\ers)_G$, $r,s \in \bN$, it follows by definition of $\cog$ that the dual $G$-action $\mu$ can be regarded as an injective $\zro$-morphism $\mu : ~ \cog \ra \mA$, such that
\[
\mu \circ \sigma_G = \rho \circ \mu \ \ , \ \ 
\mu (\sgrs) \subseteq (\rhors) \ , \ r,s \in \bN \ .
\]
\noindent Thus, the definition of dual action by Doplicher and Roberts in \cite[\S 2]{DR89A} is included in our general setting. Note that differently from the above-cited reference, our group $G$ may be noncompact: for example, consider $\mE := X \times \bC^d$, and $G$ be the group of continuous $\ud$-valued maps coinciding with the identity over some closed $W \subset X$. This will be reason of some technical complications, since no Haar measure is avalaible (Lemma \ref{lem38}, Lemma \ref{def_ra}, Lemma \ref{lem_fpa}).

Aim of the present section is to construct a universal \sC algebra $\mc$ with a commutative diagram
\begin{equation}
\label{30}
\xymatrix{
            \mA
	    \ar[r]^{i}
	 &  \mc
	 \\ \cog
	    \ar[u]^{\mu}
	    \ar@{^{(}->}[r]
	 &  \coe
	    \ar[u]_{j}
}
\end{equation}
\noindent such that
\begin{equation}
\label{eq_inn_end}
\psi \cdot i (a) = i \circ \rho (a) \cdot \psi \ ,
\end{equation}
\noindent $\psi \in j(\wE)$, $a \in \mA$. In such a case, we will say that $\mc$ is the crossed product of $\mA$ by the dual action $\mu$. Note that by commutativity of the above diagram $i(f) = j(f)$, $f \in \zro$. If $\left\{ \psi_l \right\}$ is a set of generators for $j(\wE)$, then by (\ref {eq_inn_end}) we find $i \circ \rho (a) = \sum_l \psi_l i (a) {\psi_l}^*$, so that $\rho$ becomes inner in $\mc$.

\begin{ex} {\em 
Let $X$ be a compact Hausdorff space, $\mE \ra X$ a rank $d$ vector bundle. The circle $\bT$ acts on $\mE$ by scalar multiplication, so that there is an inclusion $\bT \subset \mUE$. By (\ref{24}), we get the circle action $\bT \ra {\bf aut}_X \coe$ (see \cite[\S 3]{Pim97}); we denote by $\coe^0$ the fixed-point algebra. $\coe^0$ can be described as an inductive limit $(\mE , \mE) \stackrel{i_1}{\hra} \cdots \stackrel{i_{r-1}}{\hra} (\mE^r , \mE^r) \stackrel{i_r}{\hra} (\mE^{r+1} , \mE^{r+1}) \stackrel{i_{r+1}}{\hra} \cdots$, where $i_r(t) := t \otimes 1$. Thus, the category $\wa \bT$ has nontrivial arrows only for $r = s$, with $(\mE^r , \mE^r)_\bT = (\mE^r , \mE^r)$. Let $\sigma_\bT := \sigma |_{\coe^0}$; since $(\mE^r , \mE^r)_\bT = (\sigma_\bT^r , \sigma_\bT^r)$, the identity automorphism of $\coe^0$ defines a dual action $\mu_\mE : \wa \bT \ra \wa \sigma_\bT$. Now, the morphisms $\coe^0 \stackrel{\mu_\mE}{\ra} \coe^0 \hra \coe$ realize the property (\ref{30}), with $\mA = \coe^0$, $\mc = \coe$; moreover, since $\sigma_\bT$ becomes inner in $\coe$, we find $$\coe \simeq \coe^0 \rtimes_{\mu_\mE} \wa \bT \ .$$ Note that if $\mL \ra X$ is a line bundle, then $((\mE \otimes \mL)^r  , (\mE \otimes \mL)^r) = (\mE^r , \mE^r) = (\sigma_\bT^r , \sigma_\bT^r)$, $r \in \bN$. Thus, we also find $\mO_{\mE \otimes \mL} \simeq \coe^0 \rtimes_{\mu_{\mE \otimes \mL}} \wa \bT$. }
\end{ex}

\begin{ex} {\em 
\label{g_dual_act}
Let $\mE \ra X$ be a vector bundle over a compact Hausdorff space $X$, $G \subseteq \mUE$ a closed group. Since $(\ers)_G = (\sgrs)$, $r,s \in \bN$, a dual action $\mu : \wa G \ra \wa \sigma_G \subset {\bf end}\cog$ is defined. The inclusion $\cog \hra \coe$ induces the commutative diagram (\ref{30}), where $\mu$, $j$ are the identity map, and $i$ is the inclusion $\cog \hra \coe$. The canonical endomorphism being inner in $\coe$, the property (\ref{eq_inn_end}) is satisfied, thus $\coe = \cog \rtimes_\mu \wa G$. More generally, let $\mE , \mE' \ra X$ be vector bundles, with $G \subseteq \mUE$, $G' \subseteq \mUE'$ such that there is an isomorphism of tensor categories $\mu : \wa G \ra \wa {G'}$. Then, $\mu$ defines a dual action $\mu : \wa G \ra \wa \sigma_{G'} \subset {\bf end} \mO_{G'}$, and $\coe = \mO_{G'} \rtimes_\mu \wa G$. For example, if $\mE$, $\mE'$ have the same rank and first Chern class, then $\coe = \mO_{\mSUE'} \rtimes_\mu \wa \mSUE$ (\cite[Cor.4.18]{Vas04}). }
\end{ex}

\begin{rem} {\em 
\label{thm_bgae}
Let $\mu : \wa G \ra \wa \rho$, $G \subseteq \mUE$, be a dual action. Then, $\rho \in \tend$ has weak permutation symmetry. In fact, note that $\theta (p) \in (\mE^r , \mE^r)_G$ for every $p \in \bP_r$, $r \in \bN$ (see \S \ref{apdx}). Thus, we can define the unitary representations $\eps (p) := \mu \circ \theta (p) \in (\rho^r , \rho^r)$, $p \in \bP_r$. }
\end{rem}

\subsection{The algebraic crossed product.}

We now proceed to the construction of our crossed product, by starting at a purely *-algebraic level. Let $\aoe$ (resp. $\aog$) denote the dense $*$-subalgebra of $\coe$ (resp. $\cog$) generated by the spaces $(\ers)$ (resp. $(\ers)_G$), $r,s \in \bN$. We fix our notations with $a,a', \dots$ for elements of $\mA$, $f,f', \ldots$ for elements of $\zro$, $t,t', \ldots$ for elements of $\aoe$, $\psi, \psi', \ldots$ for elements of $(\iota,\mE)$, $y,y', \ldots$ for elements of $\aog$. 

We consider the algebraic tensor product $\mA \odot \aoe$, and the quotient $\m$ by the relation
\begin{equation}
\label{rel}
a \mu (y) \otimes t - a \otimes yt \ .
\end{equation}
\noindent By construction, $\m$ is a $\aog$-bimodule; in particular, it is a $\zro$-bimodule with coinciding left and right action. The universality implies that for each \sC algebra $\mB$ satisfying (\ref{30},\ref{eq_inn_end}) there is a unique linear map $\m \ra \mB$, defined by $\aot \mapsto i (a) j(t)$.

Our purpose is now to define a *-algebra structure on $\m$, by using the left regular representation (on the line of \cite[\S 3]{DR89A}). Let ${\bf end}(\m)$ denote the algebra of $\zro$-bimodule maps of $\m$; for each $a , \psi$, we define the following elements of ${\bf end}(\m)$:
\begin{equation}
\label{33e1}
\ovl a \ (a' \otimes_\mu t) := a a' \otimes_\mu t
\end{equation}
\begin{equation}
\label{33e2}
\ovl \psi \ (a' \otimes_\mu t) := \rho (a') \otimes_\mu \psi t \ .
\end{equation}
\noindent Note that the following relations hold:
\begin{equation}
\label{34}
\ovl \psi \ \ovl a = \ovl {\rho (a)} \ \ovl \psi \ .
\end{equation}

We now make the assumption that $G \subseteq \mSUE$, so that $(\iota , \lambda \mE) \subseteq (\iota , \mE^d)_G$. We consider a set of generators $\left\{ R_i \right\}$ of $(\iota , \lambda \mE)$, and define
\begin{equation}
\label{35}
\plihat := \sqrt d \psi_l^* R_i \in (\iota , \sigma^{d-1}) \ .
\end{equation}
\noindent Let us establish some useful relations for the $\plihat$'s. It follows by definition that $\plihat{}^* \pmjhat \in (\ii) = \zro$; furthermore, by multiplying on the left by $\psi_l$, and by summing over $l$ in (\ref {35}), we obtain
\begin{equation}
\label{36}
\sum \limits_l \psi_l \plihat = \sqrt d  R_i \ \ ;
\end{equation}
\noindent now, by definition of $\plihat$,
\[
R_i^* \psi_l = \frac 1{\sqrt d} \plihat{}^*
\]
\noindent so that, by using (\ref {eq_scp1}) we find
\[
(-1)^{d-1} \frac 1{d} \psi_l = \sum \limits_i R_i^* \sigma (R_i) \psi_l 
= \sum \limits_i R_i^* \psi_l R_i 
= \frac 1{\sqrt d} \sum \limits_i \plihat{}^* R_i \ ,
\]
\noindent hence
\begin{equation}
\label{37}
\psi_l^* = \sqrt d (-1)^{d-1} \sum \limits_i R_i^* \plihat \ \ .
\end{equation}
\noindent The above identities suggest to define (keeping in mind (\ref {33e2}), and recalling that every $R_i$ is an element of $^0 \mO_{\mSUE} \subseteq \aog$),
\begin{equation}
\label{eq_psis}
\ovl {\psi_l^*} \ (\aot) :=
\sqrt d (-1)^{d-1} \sum \limits_i \mu ( R_i^*) \rho^{d-1} (a) \otimes_\mu \plihat t \ .
\end{equation}
\noindent By the previous definition, using the fact that $\mu (R_i^*) \in (\rho^d , \iota)$, and by (\ref {37}), we obtain
\[
\ovl {\psi_l^*} ( \rho (a) \otimes_\mu t) = a \otimes_\mu {\psi_l^*} t \ .
\]

\begin{lem}
\label{lem31}
Let $\left\{ \psi_l \right\}$ be a set of generators of $\wE$, and $id \in {\bf end}(\m)$ denote the identity. The following relations hold:
\[
\ovl f \ \ovl \psi_m = \ovl \psi_m  \ovl f \ \ , \ \ 
\ovl {\psi_l^*} \ \ovl \psi_m = 
\ovl { \left\langle \psi_l , \psi_m \right\rangle } \ \ , \ \ 
\sum \limits_l \ovl \psi_l \ \ovl {\psi_l^*} = id \ \ .
\]
\end{lem}

\begin{proof} 
The first identity follows from (\ref{34}). In order to prove the second identity, we compute
\[
\begin{array}{ll}
\ovl {\psi_l^*} \cdot \ovl \psi_m ( \aot ) & = 
\ovl {\psi_l^*} ( \rho (a) \otimes_\mu \psi_m t ) = \\ & =
a \otimes_\mu (-1)^{d-1} \sqrt d \sum \limits_i R_i^* \plihat \psi_m t = \\ & =
a \otimes_\mu \psi_l^* \psi_m t \ ;
\end{array}
\]
\noindent about the third one
\[
\begin{array}{ll}
\sum \limits_l \ovl \psi_l \ovl \psi_l^* ( \aot ) & =
(-1)^{d-1} \sqrt d \sum \limits_{l,i} \ovl \psi_l 
     ( \mu (R_i^*) \rho^{d-1} (a) \otimes_\mu \plihat t ) = \\ & =
(-1)^{d-1} \sqrt d \sum \limits_{l,i} \rho \circ \mu ( R_i^* ) \rho^d (a)
     \otimes_\mu \psi_l \plihat t = \\ & =
(-1)^{d-1} d \sum \limits_i \mu \circ \sigma (R_i^*) \mu (R_i) \aot = \\ & =
\aot \ .
\end{array}
\]
\end{proof}

By (\ref{genrel}), and the previous lemma, we obtain:

\begin{cor}
\label{cor32}
The map $\psi \mapsto \ovl \psi $, $\psi \in \wE$, is a $\zro$-module map and extends to an algebra monomorphism $\aoe \hra {\bf end}(\m)$.
\end{cor}

Now, we have algebra monomorphisms from $\aoe , \mA$ into ${\bf end}(\m)$. In order to obtain the desired immersion $\m \hra {\bf end}(\m)$, we have to verify that $\ovl y = \ovl{\mu (y)}$, $y \in \cog$. For this purpose, we prove the next two lemmas.

\begin{lem}
\label{lem330}
Let $r \in \bN$. If $y \in (\iota , \mE^r)_G$, then $\ovl y = \ovl{\mu (y)}$. In particular, $\ovl{R} = \ovl{\mu(R)}$ for every $R \in (\iota , \lambda \mE) \subset (\iota , \mE^d)_G$.
\end{lem}

\begin{proof}
If $y \in (\iota , \mE^r )_G$, then $\mu (y) \in (\iota , \rho^r )$, so that
\[
\ovl y (\aot) = \rho^r (a) \otimes_\mu yt  =
\rho^r (a) \mu (y) \otimes_\mu t = \mu (y) a \otimes_\mu t =
\ovl {\mu (y)} (\aot) \ .
\]
\end{proof}

\begin{lem}
\label{lem331}
Let $\theta \in (\mE^2 , \mE^2)_G$ denote the flip (\ref{defsim}); then, $\ovl \theta = \ovl{\mu(\theta)}$.
\end{lem}

\begin{proof}
Let us compute
\[
\begin{array}{ll}
\ovl \psi_l^* \ovl \psi_m^* (\aot) & =
(-1)^{d-1} d \sum \limits_j \ovl \psi_m^* ( \mu (R_j^*) \rho^{d-1} (a) \otimes_\mu
     \pmjhat t ) = \\ & =
d \sum \limits_{i,j} \mu (R_i^*) \rho^{d-1}( \mu (R_j^*) \rho^{d-1}(a) ) \otimes_\mu
     \plihat \pmjhat t = \\ & =
d \sum \limits_{i,j} \mu (R_i)^* \rho^{d-1}( \mu (R_j)^* \rho^{2d-2}(a) ) \otimes_\mu
     \plihat \pmjhat t \ ;
\end{array}
\]
\noindent by using (\ref{36}), the identity $\rho \circ \mu (R_i) = \mu \circ \sigma_G (R_i) \in (\rho, \rho^{d+1})$, and (\ref{eq_scp1}), we obtain
\[
\begin{array}{ll}
\ovl \theta (\aot) & = 
\sum \limits_{l,m} \ovl {\psi_l \psi_m \psi_l^* \psi_m^*} (\aot) = \\ & =
d \ovl \psi_l \sum \limits_{i,j,l,m} \rho \circ \mu (R_i^*) \rho^{2d-1} (b) \otimes_\mu
       \psi_m \plihat \pmjhat t = \\ & =
d^2 \sum \limits_{i,j} \mu (\theta) \rho^2 \circ \mu (R_i^*) \rho^{d+1} \circ \mu (R_j^*)
       \rho^{2d}(a) \otimes_\mu \sigma_G (R_i) R_j t = \\ & =
d^2 \mu (\theta) \left( \sum \limits_{i,j} \rho \circ \mu ( \sigma_G (R_i^*) 
       R_i) \mu ( \sigma_G (R_j^*) R_j) \right) \aot = \\ & = 
\mu (\theta) \aot = \ovl {\mu (\theta)} (\aot) \ .
	   
\end{array}
\]
\end{proof}

About next lemma, see the "pointwise version" \cite[Lemma 2.3]{DR89A}.

\begin{lem}
\label{lem33}
If $y \in \aog$, then $\ovl y = \ovl {\mu (y)}$.
\end{lem}

\begin{proof}
We define $\mF_\mu := \left\{ y \in \aog : \ovl y = \ovl {\mu (y)} \right\}$. By the previous lemma, we have $\theta \in \mF_\mu$, $y \in \mF_\mu$, $y \in (\iota , \mE^r)_G$, $r \in \bN$. The lemma will be proven by verifying that $\mF_\mu = \aog$. By the previous corollary, the map $y \mapsto \ovl y$ is multiplicative on $\aoe \supset \aog$, thus we find that $\mF_\mu$ is a subalgebra of $\aog$. The main technical difficulty is given by the fact that {\it a priori} $\mF_\mu$ is not a $*$-subalgebra of $\aog$. As a first step, let us compute
\[
\ovl {\sigma_G (y)} = 
\sum \limits_l \ovl \psi_l \ovl y \ovl \psi_l^* = 
\ovl {\rho \circ \mu (y)} \sum \limits_l \ovl \psi_l \ovl \psi_l^* =
\ovl {\rho \circ \mu (y)} = \ovl {\mu \circ \sigma_G (y)} \ ;
\]
\noindent we proved in this way that $\mF_\mu$ is $\sigma_G$-stable. 

We want to prove that $R^* \in \mF_\mu$, $R \in (\iota , \lambda \mE)$. For this purpose, we note that $P_{\theta,d}$, defined as in (\ref{eq_P}), belongs to the $\sigma_G$-stable algebra generated by $\theta$: since we proved that $\mF_\mu$ is a $\sigma_G$-stable algebra and contains $\theta$, we find $P_{\theta,d} \in \mF_\mu$. Let $\left\{ R_i \right\}$ be a set of generators of $(\iota , \lambda \mE)$; then, for every index $i$ there is $f_i \in C(\spzro)$ such that $R_i R_i^* = f_i P_{\theta,d}$: since $\zro \subset \mF_\mu$, we find $R_i R_i^* \in \mF_\mu$. In order to prove that each $R_i^*$ belongs to $\mF_\mu$, note that $R_i R_i^* R_i = \lambda_i^2 R_i$ (where $\lambda_i^2 := R_i^* R_i \in \zro$), and compute
\[
\begin{array}{ll}
\lambda_i^2 \ovl {R_i^*} = \ovl {R_i^* R_i R_i^*} = \ovl R_i^* 
         \ovl { \mu (R_i R_i^*) } = \ovl R_i^* \ovl {\mu (R_i)}
		 \ovl {\mu (R_i^*)} =
         \ovl R_i^* \ovl R_i \ovl {\mu (R_i^*)} = \lambda_i^2 \ovl {\mu (R_i^*)} \ .
\end{array}
\]
\noindent Since each $R_i$ is a section of $\lambda \mE$ having the same support as $\lambda_i$, we conclude that $\ovl {R_i^*} = \ovl {\mu (R_i^*)}$, i.e. $R_i^* \in \mF_\mu$. 


Let finally $y \in (\ers)_G$. Given a set $\left\{ \psi_l \right\}$ of generators for $\wE$, we consider the tensors $\psi_L \in (\iota , \sigma^r)$ (see (\ref{def_psiL})), so that $y$ is a linear combination with coefficients in $\zro$ of terms of the type $\psi_L \psi_M^* $; by using (\ref{37}), $R_i^* \in (\sigma_G^d , \iota)$, and $\plihat \in (\iota,\sigma^{d-1})$, we obtain
\[
\begin{array}{ll}
\psi_L \psi_M^* = \sum \limits_J \sigma_G^r ( R_{j_1}^* \sigma_G^{d-1} (R_{j_2}^*) \ \cdots \ 
       \sigma_G^{(s-1)(d-1)} (R_{j_s}^*) ) \ \psi_L {\wa \psi}_{m_1}^{j_1} \ \cdots \
	   {\wa \psi}_{m_s}^{j_s} \ .
\end{array}
\]
\noindent Let now $U \subseteq \spzro$ be a closed set trivializing $\lambda \mE$, $R_U$ the generator of the module of continuous sections of $\lambda \mE |_U$ (so that, $R_U$ is an isometry in $\coe |_U$); then, the restriction $y_U \in \coe |_U$ of $y$ over $U$ is given by $y_U = T_U y'_U$, where
\[
T_U := 
\sigma_G^r ( \ R_U^* \sigma_G^{d-1} (R_U^*) \ \cdots \ \sigma_G^{(s-1)(d-1)} (R_U^*) \ ) 
\ ,
\]
\noindent and
\[
y'_U :=
\sum \limits_{LM} f^U_{LM} \psi_L {\wa \psi}_{m_1}^U \ \cdots \ {\wa \psi}_{m_s}^U \ ;
\]
\noindent here $f^U_{LM}$ are suitable elements of $C(U)$, and each ${\wa \psi}_{m_k}^U$ is defined as in (\ref{35}) by replacing $R_i$ with $R_U$. Since $R_U^* \in \mF_\mu |_U$, we obtain that $T_U$ is an isometry in $\mF_\mu |_U$; moreover, $y'_U \in (\iota , \mE^{r + s (d-1)}) |_U$. Now, $T_U y'_U = y_U = {\wa g} (y_U) = {\wa g} (T_U y'_U) = T_U {\wa g} (y'_U)$; thus, by multiplying on the left by $T_U^*$, we obtain ${\wa g} (y'_U) = y'_U$. Thus, we conclude that $y'_U \in  (\iota,\mE^{r+s(d-1)})_G |_U \subset  \mF_\mu |_U$, and $y_U \in \mF_\mu |_U$. By applying the above argument for a closed trivializing cover of $\spzro$, we obtain $y \in \mF_\mu$.
\end{proof}

By the previous lemma, the maps (\ref{33e1},\ref{33e2}) extend to an injective $\zro$-bimodule map from $\m$ into ${\bf end}(\m)$:
\begin{equation}
\label{mendm}
\ovl {\aot} := \ovl a \ \ovl t  \ ,
\end{equation}
\noindent (in fact, (\ref{mendm}) factorizes through elements of $\aog$); the injectivity of (\ref{mendm}) follows from the equality $\ovl {\aot} (1 \otimes_\mu 1) = \aot$. We now can regard at $\m$ as a {\em subalgebra} of ${\bf end}(\m)$; moreover, the map 
\[
\ovl {\aot} \mapsto {\ovl {\aot}}^*  :=  \ovl {1 \otimes_\mu t^*} \cdot \ovl {a^* \otimes_\mu 1}
\] 
\noindent defines an involution on $\m$, extending the involutions defined on $\mA$, $\aoe$, as can be proven exactly in the same way used in \cite[Lemma 2.5]{DR89A}.

In order to economize in notations, in the sequel we will identify $\m$ with the image in ${\bf end}(\m)$, so that the overline symbol is dropped. We denote by $i(a) := a \otimes_\mu 1$, $j(t) := 1 \otimes_\mu t$ the canonical immersions $\mA,\aoe$ $\hra$ $\m$. The following lemma is the solution to our initial problem at the *-algebraic level.

\begin{lem}
Let $\rho$ be a unital endomorphism of a \sC algebra $\mA$, $\mE \ra \spzro$ a vector bundle, $\mu : \wa G \ra \wa \rho$ a dual action, $G \subseteq \mSUE$. Then
\begin{enumerate}

\item for each commutative diagram of *-algebra morphisms
\begin{center}
$\xymatrix{
   \mA
   \ar[r]^{\pi}
&  \mB
\\ \aog
   \ar@{^{(}->}[r]
   \ar[u]^{\mu}
&  \aoe
   \ar[u]_{\phi}
}$ \\
\end{center}
such that $\phi (f) = \pi (f)$, $\phi (\psi) \pi (a) = \pi \circ \rho (a) \phi (\psi)$, there is a unique *-homomorphism $\chi : \m \ra \mB$ such that $\chi (\aot) = \pi (a) \phi (t)$;

\item there is an action by *-automorphisms $\alpha : G \ra {\bf aut}(\m)$, $$\alpha_g (\aot) := a \otimes_\mu \wa g(t) \ .$$

\end{enumerate}
\end{lem}

\begin{proof}
We prove the first point. By universality, the application $\chi$ exists as a linear map. Thus, we have to prove that $\chi$ is a *-homomorphism, and this can be done exactly as in \cite[Thm.2.6]{DR89A}. About the second point, note that by Lemma \ref{lem33} the identity $i \circ \mu (y) = j(y)$ holds for every $y$; thus, we can construct a commutative diagram as in the first point, with $\mB = \m$, $i = \pi$, $j \circ \wa g = \phi$. Furthermore, $j \circ \wa g (f) = i (f)$ for $f \in \zro$, and
\[ 
i \circ \wa g (\psi) \cdot i (a) = 1 \otimes_\mu \wa g(\psi) \cdot a \otimes_\mu 1
    = \rho (a) \otimes_\mu 1 \cdot 1 \otimes_\mu \wa g (\psi)
    = i \circ \rho (a) \cdot j \circ \wa g (\psi) \ ;
\]
\noindent so that, we find that the property (\ref{eq_inn_end}) holds for elements of $i(\mA)$ and $j \circ \wa g (\aoe)$. Thus, by universality there is a *-endomorphism $\alpha_g \in {\bf end}(\m)$. Now, $\alpha_g$ has inverse $\alpha_{g^{*}}$, and it is obvious that $\alpha_{g{g'}} = \alpha_g \alpha_{{g'}}$.
\end{proof}

\subsection{The $\zro$-crossed product.}

Our aim is now to extend the previous lemma in such a way to get (\ref{30},\ref{eq_inn_end}) in the setting of \sC algebra morphisms. It is clear that the solution will be given by the closure of the *-algebra $\m$ w.r.t. some suitable \sC norm. For this purpose, we construct a family of \sC seminorms for $\m$, indexed by elements of the spectrum of $\zro$.

Let $x \in \spzro$, and $\ker x$ be the corresponding ideal of $\zro$. We consider the ideal $\mA \ker x \subset \mA$; furthermore, we introduce the notation $\mA_x$ to indicate the quotient of $\mA$ by $\mA \ker x$, and denote by $\pi_x : \mA \ra \mA_x$ the natural surjection. By \cite[Prop.2.8]{Bla96}, we have $\left\| a \right\| = \sup_x \left\| \pi_x (a) \right\|$ for each $a \in \mA$. Let now $G^x \subseteq \sud$ be defined as by (\ref{eq_def_G^x}). We denote by $x_* : \coe \ra \mO_d$ the evaluation at $x$ of $\coe$ as a continuous bundle of \sC algebras. We want to construct a morphism $\mu_x : \mO_{G^x} \ra \mA_x$ such that $\pi_x \circ \mu = \mu_x \circ x_*$. In the same way, we want to define an endomorphism $\rho_x : \mA_x \ra \mA_x$ such that $\pi_x \circ \rho = \rho_x \circ \pi_x$. In order to define $\mu_x$, we will make use of the following lemma:

\begin{lem}
Let $( \mF , \mF \ra \mF_x)$ be a continuous bundle of \sC algebras over a compact space $X$, $U \subseteq X$ a closed set. If $y \in \mF$ is a vector field such that the restriction $y_U$ vanishes over $U$, then there exists a continuous function $f \in C(X)$ such that $f |_U = 0$, and a vector field $y' \in \mF$ with $y = fy'$.
\end{lem}

\begin{proof}
It suffices to regard at the \sC algebra $\mF_U := \left\{ y \in \mF : y_U = 0  \right\}$ as a non-degenerate Banach $C_0(X - U)$-bimodule, and apply \cite[Prop.1.8]{Bla96}.
\end{proof}

\noindent We now give the following definitions:
\begin{equation}
\mu_x \circ x_* (y) := \pi_x \circ \mu (y)
\end{equation}
\begin{equation}
\label{def_ro}
\rho_x \circ \pi_x(a) := \pi_x \circ \rho (a) \ .
\end{equation}
\noindent Since $\mu$ is a $\zro$-module map, it follows from the previous lemma that if $y_x := x_* (y) = 0$, then $\pi_x (\mu (y)) = f(x) \cdot \pi_x (\mu (y)) = 0$, so that $\mu_x$ is well defined. About $\rho_x$, note that if $\pi_x(a) = 0$, then $a$ is norm limit of terms of the type $\sum_k f_k a_k$, where $f_k \in \ker x \subset \zro$, $a_k \in \mA$. Thus, $\pi_x (\rho (a)) = 0$ (in fact, the restriction of $\pi_x$ to $\zro$ is the evaluation on $x$), and $\rho_x$ is well defined. We can now prove the following

\begin{prop}
\label{point_dual_action}
For each $x \in \spzro$, the morphism $\mu_x : (\cog)_x \ra  \mA_x$ defines a dual $G^x$-action $\wa {G^x} \ra \tend_x$, where $G^x \subseteq \sud$ is the closed group defined in (\ref{eq_def_G^x}).
\end{prop}

\begin{proof}
We have constructed for every $x \in \spzro$ a morphism $\mu_x : (\cog)_x \simeq \mO_{G^x} \ra \mA_x$. First of all, we have to verify that $\rho_x \circ \mu_x = \mu_x \circ \sigma_x$, where $\sigma_x$ is the canonical endomorphism of $\mO_{G^x}$. If $y_x \in \mO_{G^x}$, we find $\rho_x \circ \mu_x (y_x) = \rho_x \circ \pi_x \circ \mu (y)  = \pi_x \circ \rho \circ \mu (y) = \pi_x \circ \mu \circ \sigma (y) = \mu_x \circ \pi_x \circ \sigma (y) = \mu_x \circ \sigma_x (y_x)$. We now prove that $\mu_x (\sigma_x^r,\sigma_x^s) \subseteq (\rho_x^r,\rho_x^s)$, $r,s \in \bN$: if $y_x \in (\sigma_x^r,\sigma_x^s)$, then for every $\pi_x (a) \in \mA_x$ we find $\mu_x (y_x) \cdot \rho_x^r \circ \pi_x (a) = \pi_x \circ \mu (y) \cdot \pi_x \circ \rho^r (a) = \pi_x (\mu (y) \cdot \rho^r (a) ) = \pi_x ( \rho^s(a) \cdot \mu(y) ) = \rho_x^s \circ \pi_x (a) \cdot \mu_x (y_x)$, and the proposition is proven.
\end{proof}

We can now apply the results of \cite[\S 3]{DR89A}, and define $\mo$ as the closure w.r.t. the maximal \sC seminorm of the *-algebra $\mA_x \odot_{\mu_x} {^0\mO_d}$, obtained by the dual $G^x$-action defined in the previous proposition. By \cite[Thm.3.2]{DR89A}, $\mo$ is the solution of our universal problem in the case of the 'flat' dual action $\mu_x : \wa {G^x} \ra \wa \rho_x$; moreover, there is an action $\alpha_x : G^x \ra {\bf aut} (\mo)$, $x \in \spzro$, such that 
\begin{equation}
\label{eq_alphax}
\alpha_{x,u} (a_x \otimes_{\mu_x} t_x) = a_x \otimes_{\mu_x} \wa u (t_x) \ \ ,
\end{equation}
\noindent $u \in G^x$, $a_x \in \mA_x$, $t_x \in \mO_d$.

\begin{lem}
\label{seminorms}
For each $x \in \spzro$, there is a *-homomorphism ${\wa x} : \m \ra \mo$ extending the evaluation morphisms $\pi_x$, $x_*$.
\end{lem}

\begin{proof}
We define $\wa x (\aot) := \pi_x(a) \otimes_{\mu_x} x_*(t)$. In order to verify that $\wa x$ is well defined, we have to prove that $\wa x ( a \mu(y) \otimes_\mu t ) = \wa x ( a \otimes_\mu yt )$. This immediately follows from the definition of $\mu_x$.
\end{proof}

We introduce on $\m$ the \sC seminorm $\left\| \cdot  \right\| := \sup_x  \left\| \wa x ( \cdot ) \right\|$, and denote by $\mc$ the closure of $\m$ w.r.t. $\left\| \cdot  \right\|$. By construction $ \left\| a \right\| = \left\| a \otimes_\mu 1 \right\| $, and $\left\| t \right\| = \left\| 1 \otimes_\mu t \right\|$ for $a \in \mA$, $t \in \coe$. Thus, there are unital monomorphisms $i : \mA \hra \mc$ and $j : \coe \hra \mc$; {\em in the sequel, we will identify elements of $\coe$, $\mA$ with their images in $\mc$}. Note that in general $\mc$ {\em is not} a continuous bundle of \sC algebras over $\spzro$, but simply a $\zro$-algebra. For every $x \in \spzro$, there are epimorphisms
\begin{equation}
\label{eq_epix}
\eta_x : \mc \ra \mo \ \ .
\end{equation}

%
%
%
%
%

\begin{thm}
\label{thm38}
Let $\mA$ be a unital \sC algebra, $\mE \ra \spzro$ a vector bundle, $G \subseteq \mSUE$ a closed group with a dual action $\mu : \wa G \ra \wa \rho \subset \tend$. Then, $\mc$ is the unique crossed product of $\mA$ by $\mu$ satisfying the universal properties (\ref{30}, \ref{eq_inn_end}), and carrying a strongly continuous action $\alpha : G \ra {\bf aut}_{\spzro}(\mc)$, $\alpha_g (\aot) = a \otimes_\mu \wa g (t)$.
\end{thm}

\begin{proof}
We prove the unicity of the \sC norm. Let $\mB$ be the completition of $\m$ w.r.t. a given \sC norm, extending those defined on $\mA$, $\coe$; then, $\mB$ is a $\zro$-algebra (i.e., an upper semicontinuous bundle over $\spzro$), with fibre epimorphisms $\phi_x : \mB \ra \mB_x$, $x \in \spzro$. Since $\phi_x$ extends the evaluation epimorphisms $\pi_x : \mA \ra \mA_x$, $x_* : \coe \ra \mO_d$, every $\mB_x$ is the closure of $\mA_x \odot_{\mu_x} {^0 \mO_d}$ w.r.t. some \sC norm. Let $g \in G$; since $\alpha_g \in {\bf aut}_{\spzro} \mB$, a field of automorphisms $\alpha_{x,g} \in {\bf aut} \mB_x$ is defined (\cite[\S 4]{Nil96}), with $\alpha_{x,g} \circ \phi_x = \phi_x \circ \alpha_g$. Let $g(x) \in G^x \subseteq \sud$ denote the evaluation of $g \in G$ over $x \in \spzro$. By definition, 
\begin{equation}
\label{eq_ax}
\alpha_{x,g} ( \pi_x(a) \otimes_{\mu_x} x_*(t) ) = 
\phi_x (a \otimes_\mu \wa g (t)) =
\pi_x(a) \otimes_{\mu_x} \wa{g(x)} (x_*(t))  \ ,
\end{equation}
\noindent $a \in \mA$, $t \in \aoe$. Now, by \cite[Thm.3.2]{DR89A}, the unique \sC norm on $\mA_x \odot_{\mu_x} \mO_d$ such that the r.h.s. of (\ref{eq_ax}) defines an (isometric) automorphism is the maximal \sC norm, thus $\mB_x = \mo$. By a general property of upper semicontinuous bundles, we must have $\left\| b \right\| = \sup_x \left\| \phi_x (b) \right\|$, $b \in \mB$; since every $\left\| \phi_x (\cdot) \right\|$, $x \in \spzro$, is the maximal seminorm over $\mA_x \odot_{\mu_x} {^0 \mO_d}$, we conclude that $\mB = \mc$. Note that with the notation of (\ref{eq_alphax},\ref{eq_epix}), we also proved 
\begin{equation}
\label{eq_intact}
\eta_x \circ \alpha_g = \alpha_{x , g(x) } \circ \eta_x
\ \ , \ \
x \in \spzro \ .
\end{equation}
\end{proof}

It is clear that $\mA$ is contained in the fixed-point algebra of $\mc$ w.r.t. the $G$-action. Let us now suppose that $G$ is a compact group w.r.t. the natural topology induced as a subgroup of $\mUE$. Then, the Haar measure on $G$ defines an $\mA$-valued invariant mean on $\m$: $$m (\aot) := a \ \mu \left( \int_G \wa g (t) \ dg \right) \ .$$ The argument of \cite[Lemma 2.8, Thm.3.2]{DR89A} now applies, and we conclude that $\mA = (\mc)^\alpha$. In the general case, in order to prove that $\mA$ is the fixed-point algebra w.r.t. the $G$-action some analogue of the Haar mean is needed; we will return on this point in Lemma \ref{lem38}.

\begin{ex} {\em
Let $G$ be a compact group, $\mE \ra X$ a $G$-vector bundle with trivial action on $X$ in the sense of \cite{Ati}. Then, we may regard at $G$ as a closed (compact) subgroup of $\mUE$, and the previous argument applies.
}
\end{ex}

\begin{rem} {\em 
\label{rem_ce}
Let $\left\{ \psi_l \right\}$ be a set of generators for $\wE$. We denote by $\sigma_\mE \in {\bf end}_{\spzro} (\mc)$, $\sigma_\mE (b) := \sum_l \psi_l b \psi_l^*$, $b \in \mc$, the inner endomorphism induced by $\wE$. By definition, $\wE \subseteq (\iota , \sigma_\mE)$. By (\ref{eq_inn_end}), $\sigma_\mE |_\mA = \rho$. Also note that $\sigma_\mE |_{\coe} = \sigma$, where $\sigma \in {\bf end}_X \coe$ is defined by (\ref{def_ce}). }
\end{rem}

\begin{rem} {\em 
Since $\mc$ satisfies (\ref{eq_inn_end}), we find that $\mc$ is a quotient of the crossed product $\cpen$ in the sense of \cite[Prop.3.4]{Vas04}. }
\end{rem}

Crossed products by dual actions with minimal relative commutant of $\mA$ will be of particular interest for our purposes. We denote by ${\bf aut}_\mA (\mc , \sigma_\mE)$ the group of automorphisms on $\mc$ leaving $\mA$ pointwise fixed, and commuting with $\sigma_\mE$. We recall that if $\mA$ is a $C(X)$-algebra, then a {\em nc-pullback} of $\mE \ra X$ is a Hilbert $\mA$-bimodule $\mM$ which is isomorphic {\em as a right Hilbert $\mA$-module} to the tensor product $\wE \otimes_{C(X)} \mA$ with coefficients in $C(X)$ (\cite[Def.5.4]{Vas04}).

\begin{prop}
\label{cp_mcr}
Let $\rho$ be a unital endomorphism of a \sC algebra $\mA$, $\mE \ra \spzro$ a vector bundle with a dual action $\mu : \wa G \ra \wa \rho$, $G \subseteq \mSUE$. Suppose $\mA' \cap (\mc) = \mA' \cap \mA =: \mZ$. Then, the Banach space
\begin{equation}
\label{def_bm}
\mM := \left\{ 
\varphi \in \mc \ : \ \varphi a = \rho (a) \varphi \ , \ a \in \mA  
\right\}
\end{equation}
\noindent is a $G$-Hilbert $\mZ$-bimodule in $\mc$, and a nc-pullback of $\mE$. Finally, there is a group isomorphism $G \simeq {\bf aut}_\mA (\mc , \sigma_\mE)$.
\end{prop}

\begin{proof}
It is clear that $\mM$ is closed under left and right multiplications by elements of $\mZ$; moreover, $\varphi^* \varphi' \in \mA' \cap (\mc)$, $\varphi , \varphi' \in \mM$, thus $\varphi^* \varphi' \in \mZ$ and $\mM$ is a Hilbert $\mZ$-bimodule in $\mc$. In order to prove that $\mM$ is a nc-pullback of $\mE$, it suffices to verify that $\mM = \wE \cdot \mZ := {\mathrm{span}} \left\{ \psi z , \psi \in \wE , z \in \mZ \right\}$: now, it is clear that $\wE \cdot \mZ \subseteq \mM$; viceversa, let $\left\{ \psi_l \right\}$ be a set of generators for $\wE$; then $\psi_l^* \varphi \in \mA' \cap (\mc) = \mZ$, thus $\varphi = \sum_l \psi_l (\psi_l^* \varphi) \in \wE \cdot \mZ$, and $\mM$ is a nc-pullback. Now, $z = \alpha_g (z)$, $z \in \mZ$, $g \in G$; thus, $\alpha_g (z \varphi z') = z \alpha_g (\varphi) z'$, $z' \in \mZ$, and $\alpha_g (\varphi^* \varphi') = \varphi^* \varphi'$, $\varphi , \varphi' \in \mM$ (in fact, $\varphi^* \varphi' \in \mZ$); moreover, $\alpha_g (\varphi) = \sum_l \psi_l (\psi_l^* \alpha_g(\varphi))$, with $\psi_l^* \alpha_g(\varphi) a = \psi_l^* \alpha_g (\rho(a) \varphi) = a \psi_l^* \alpha_g(\varphi)$, $a \in \mA$; thus $\psi_l^* \alpha_g(\varphi) \in \mZ$, and $\alpha_g(\varphi) \in \mM$. We conclude that $G$ acts on $\mM$ by unitary $\mZ$-bimodule maps. Finally, by construction, if $g \in G$ then $\alpha_g \in {\bf aut}_\mA (\mc , \sigma_\mE)$ (recall Rem.\ref{rem_ce}). Viceversa, let $\alpha \in {\bf aut}_\mA (\mc , \sigma_\mE)$, $\psi , \psi' \in \wE$; then, a direct computation shows that $\psi^* \alpha (\psi') a = a \psi^* \alpha (\psi')$, $a \in \mA$. Thus, $\psi^* \alpha (\psi') \in \mZ$. Moreover, since $\sigma \in {\bf end}_X \coe$ has permutation symmetry (see \cite[(4.16)]{Vas04}), we find $\sigma (\psi) = \theta \psi$, where $\theta \in \cog$ is defined by (\ref{defsim}). By using the relation $\eps = \mu (\theta)$ (see Rem.\ref{thm_bgae}), we find $\sigma_\mE (\psi) = \eps \psi$, thus
\[
\rho ( \psi^* \alpha (\psi') ) =
\sigma_\mE (\psi)^*  \rho \circ \alpha (\psi') =
\psi^* \eps^* \cdot \alpha \circ \sigma_\mE (\psi') =
\psi^* \eps^* \alpha ( \eps \psi' ) =
\psi \eps^* \eps \alpha (\psi') \ .
\]
\noindent  We conclude that $\psi^* \alpha (\psi') \in \zro$; thus, $$g := \sum_{lm} \psi_l \psi_m^* [\psi_l^* \alpha (\psi_m)]$$ belongs to $\mUE$, and $\alpha = \alpha_g$ (in fact, $\alpha (\aot) = a \otimes_\mu \wa g (t)$, $a \in \mA$, $t \in \aoe$). Since $\alpha_g |_{\cog}$ is the identity (recall $\cog \subset \mA$), we conclude $g \in G$ (in fact, $G$ is the stabilizer of $\cog$ in $\coe$, see \S \ref{apdx}).
\end{proof}

\begin{rem} {\em 
With the notation of the previous proposition, $G$ acts on $\mM$ by a tensor action in the sense of \cite[Def.5.3]{Vas04}. }
\end{rem}

\begin{rem} {\em 
An interesting point is the computation of the $K$-theory of $\mc$ in terms of the $K$-theory of $\mA$, $\cog$, $\coe$. By functoriality, the commutative diagram (\ref{30}) induces a sequence at level of corresponding $K$-groups
\[
\xymatrix{
            K_*(\mA)
	    \ar[r]^{i_*}
	 &  K_*(\mc)
	 \\ K_*(\cog)
	    \ar[u]^{\mu_*}
	    \ar[r]
	 &  K_*(\coe)
	    \ar[u]_{j_*}
}
\]
\noindent with $* = 0,1$. For this purpose, the $K$-theory of $\coe$ has been computed in the case in which $X$ is a $n$-sphere or a Riemann surface (\cite{Vas02}). The $K$-theory of $\cog$ should depend only on the tensor \sC category $\wa G$, and could be related to the $K$-theory of $\wa G$ in the sense of \cite{Mit01}. For $X = \left\{ x \right\}$, the $K$-theory of $\cog$ has been (implicitly) computed in \cite{PR96}. }
\end{rem}

\section{Special Endomorphisms.}
\label{spec_end}

\subsection{Basic properties.}
We start the present section by giving some basic properties about endomorphisms carrying a permutation symmetry, generalizing the analogues in \cite[\S 4]{DR87}. For our purposes, it will be convenient to start from the notion of {\em weak permutation symmetry} (Def.\ref{def11}) for an endomorphism $\rho$ of a unital \sC algebra $\mA$. In fact, in the nontrivial centre case the property (\ref{ps3}) is too restrictive, since it implies $\rho (z) = z$, $z \in \mZ := \mA \cap \mA'$.

\begin{defn}
\label{def_strict_int}
Let $\rho$ be a unital endomorphism on a \sC algebra $\mA$ carrying a weak permutation symmetry $p \mapsto \eps (p)$, $p \in \bP_\infty$. Are said to be {\bf symmetry intertwiners} the elements of $(\rhors)$ for which (\ref{ps3}) holds:
\[
( \rhors )_\eps := 
\left\{ t \in (\rhors) : \eps(s,1) t = \rho(t) \eps (r,1) \right\} \ .
\]
\end{defn}

\noindent We denote by $\wa \rho_\eps$ the \sC subcategory of $\wa \rho$ with objects the tensor powers $\rho^r$ and arrows $( \rhors )_\eps$.

\begin{lem}
\label{strict_cat}
The set of symmetry intertwiners is closed for multiplication and $\rho$-stable. So that $\wa \rho_\eps$ is a symmetric tensor \sC category.
\end{lem}

\begin{proof}
As first, we prove that if $t \in ( \rhors )_\eps , t' \in ( \rho^s , \rho^h )_\eps$, then $t't \in ( \rho^r , \rho^h )_\eps$. This fact is obvious by (\ref{ps3}), which implies 
\[
\rho(t't) = \eps(h,1) t' \eps(1,s) \eps(s,1) t \eps(1,r) \ .
\]
\noindent We now prove that for each $k \in \bN$ the inclusion $( \rhors )_\eps \subset ( \rho^{r+k} , \rho^{s+k} )_\eps$ holds. For this purpose, note that $(1,s+1) (s,1) = \bS^s(1,1)$, so $\eps(1,s+1) \eps(s,1) = \rho^s(\eps)$, and $\eps(s+1,1) \rho^s(\eps) t = \rho(t) \eps(r+1,1) \rho^r(\eps)$ for $t \in {(\rhors)_\eps}$. Since $\rho^s(\eps) t = t \rho^r(\eps)$ and $\rho^r(\eps^2) = 1$, we find that $t \in ( \rho^{r+1} , \rho^{s+1})_\eps$. The above computations imply that the set of symmetry intertwiners is closed for multiplication. In order to prove that $\wa \rho_\eps$ is $\rho$-stable, it suffices now to verify that $\eps(r,1) \in ( \rho^{r+1},\rho^{r+1} )_\eps$; but this fact is obviously true, since the relation $(r+1,1)(r,1)(1,r+1) = \bS (r,1)$ implies that $\eps(r,1)$ satisfies (\ref{ps3}). The symmetry follows from the fact that (\ref{ps3}) implies $t \rho^{r'}(t') = \eps (r',s') t' \rho^r(t) \eps(r,s)$, $t \in ( \rhors )_\eps$, $t \in ( \rho^{r'} , \rho^{s'} )_\eps$.
\end{proof}

We denote by $\soro$ the \sC algebra generated by the $( \rhors )_\eps$'s. It follows from the previous lemma that $\soro$ is $\rho$-stable. Note that $( \ii )_\eps = \zro$, so that $\soro$ is a $\zro$-algebra. Also note that in general $\wa \rho_\eps$ is not a DR-category (\cite{DR89}): in fact, in general $\wa \rho_\eps$ does not have special objects in the sense of Doplicher-Roberts \cite[\S 3]{DR89}; moreover, in general $(\ii)_\eps \neq \bC$ (for a so-called DR-category, the condition $(\ii) = \bC$ is required).

\begin{rem} {\em 
\label{rem_si_da}
Symmetry intertwiners appear when there is a dual action $\mu : \wa G \ra \wa \rho$, $G \subseteq \mUE$. In this case, $\rho$ has weak permutation symmetry (Rem.\ref{thm_bgae}), and $\mu (\ers)_G \subset ( \rhors )_\eps$ for each $r,s$. In fact, $\rho$ has permutation symmetry $\eps (p) := \mu \circ \theta (p)$, $p \in \bP_\infty$, so that $\eps (s,1) \mu (t) = \rho \circ \mu(t) \eps (r,1)$, $t \in (\ers)_G$. In particular, every endomorphism carrying weak permutation symmetry has symmetry intertwiners, thanks to the naturally induced dual $\ud$-action (see Lemma \ref{lem_da} in the present section). }
\end{rem}

We can now give the notion of permutation symmetry in our general setting, naturally coinciding with Def.\ref{def11} when formulated in the trivial centre case. In the general case, the above remarks imply that Def.\ref{def11} is more restrictive w.r.t. the following

\begin{defn}[quasi-symmetry]
\label{def_gps}
A unital endomorphism $\rho$ of a \sC algebra $\mA$ with centre $\mZ$ has permutation quasi-symmetry if there is a unitary representation $p \mapsto \eps (p)$ of the group $\bP_\infty$ of finite permutations of $\bN$ in $\mA$, such that:
\begin{equation}\label{gps1} \eps (\bS p) = \rho \circ \eps (p) \end{equation}
\begin{equation}\label{gps2} \eps := \eps (1,1) \in (\rho^2 , \rho^2) \end{equation}
\begin{equation}\label{gps3} (\rhors) = \rho^s (\mZ) \cdot (\rhors)_\eps = (\rhors)_\eps \cdot \rho^r (\mZ) \ . \end{equation}
\end{defn}

\noindent Thus, the notion of weak permutation symmetry is not affected by Def.\ref{def_gps}, and remains the one given in Def.\ref{def11}. The property (\ref{gps3}) has to be intended in the sense that $(\rho^r,\rho^s)$ is generated as a Banach space by elements belonging to $\rho^s (\mZ) \cdot (\rhors)_\eps = (\rhors)_\eps \cdot \rho^r (\mZ)$. We say that $\rho$ has {\em trivial} permutation symmetry when $\eps(p) \equiv 1$ for every $p \in \bP_\infty$. In this case $(\rhors)_\eps = \left\{ t \in (\rhors) : \rho (t) = t \right\}$ for every $r,s \in \bN$.

\bigskip

\begin{ex} {\em 
\label{ps_bp}
Let $\mE \ra X$ be a vector bundle, $\mZ$ an abelian, unital $C(X)$-algebra, $\mM \simeq \wE \otimes_{C(X)} \mZ$ a nc-pullback, $G \subseteq \mUE$ a closed group acting on $\mM$ by a tensor action in the sense of \cite[Def.5.4]{Vas04}. Then, with the notation of \cite[\S 5.1]{Vas04}, the canonical endomorphism $\tau_G$ of $\cmg$ has weak permutation symmetry in the sense of Def.\ref{def_gps}. If $G$ satisfies a suitable local triviality, then $\tau_G$ has permutation quasi-symmetry (\cite[Cor.5.8]{Vas04}). }
\end{ex}

We now establish some basic properties of endomorphisms carrying (weak) permutation symmetry. For this purpose, recall that if $\mE \ra X$ is a rank $d$ vector bundle, then there is an isomorphism $\mO_{\mUE} \simeq C(X) \otimes \mO_{\ud}$ (\cite[Prop.4.15]{Vas04}), where $\mO_{\ud} \subset \mO_d$ is the fixed-point algebra w.r.t. the canonical $\ud$-action; in the sequel, we will identity $\mO_{\mUE}$ with $C(X) \otimes \mO_{\ud}$. Also recall that $P_{\theta,d} \in (\mE^d , \mE^d )_{\mUE} \subset \mO_{\mUE}$ (\cite[\S 4.2]{Vas04}). We denote by $\delta \in {\bf end}_X \mO_{\mUE}$ the shift
\[
\delta (t) := P_{\theta,d} \otimes t \ \ , \ \ 
t \in (\ers)_{\mUE} \ , \ r,s \in \bN \ \ .
\]

\begin{lem}
\label{cp_mu}
Let $\mA$ be a unital $C(X)$-algebra, $\mu_0 : C(X) \otimes \mO_{\ud} \ra \mA$ a $C(X)$-morphism. Suppose there is a finitely generated Hilbert $C(X)$-bimodule $\mR \subset \mA$, such that
\[
\left\{ 
\begin{array}{l}
\mR \mR^* := 
{\mathrm{closed \ span}} \left\{ R'R^* , R,R' \in \mR \right\} =
C(X)  \mu_0 (P_{\theta,d})  \\ 
R \mu_0 (y) = \mu_0 \circ \delta(y) \ R  \ \  , \ \  
y \in C(X) \otimes \mO_{\ud} \ , \ R \in \mR \ .
\end{array}
\right. 
\]
\noindent Then, for every vector bundle $\mE \ra X$ with a bimodule isomorphism $\beta : (\iota , \lambda \mE) \stackrel{\simeq}{\ra} \mR$, $\mu_0$ uniquely extends to a morphism $\mu : \mO_{\mSUE} \ra \mA$, by defining $\mu (y) := \beta (y)$, $y \in (\iota , \lambda \mE)$.
\end{lem}

\begin{proof}
By \cite[\S 1.6]{Kas88}, there is a unital $C(X)$-monomorphism $\mA \hra L(\mM)$, where $\mM$ is a Hilbert ${C(X)}''$-module. Thus, for every vector bundle $\mE$ with $(\iota , \lambda \mE) \simeq \mR$, we can regard at $\mu_0$ as a covariant representation of $(\mO_{\mUE} , \delta)$ with rank $\lambda \mE$ (in the sense of \cite[Def.3.2]{Vas04}). Since $\mO_{\mSUE}$ is the crossed product of $\mO_{\mUE}$ by $\delta$ (\cite[Prop.4.17]{Vas04}), there is a $C(X)$-morphism $\mu : \mO_\mSUE \ra L(\mM)$ having image generated by $\mu_0(C(X) \otimes \mO_{\ud})$ and $\mR$, so that $\mu (\mO_\mSUE) \subseteq \mA$.
\end{proof}

By definition, every endomorphism $\rho \in \tend$ carrying weak permutation symmetry $p \mapsto \eps(p)$ induces a morphism $C^*(\bP_\infty) \ra \mA$. We now consider the finite set of characteristic functions $\left\{ \chi_i \right\} \subset \zro$ corresponding to the decomposition of $\spzro$ into connected components, and decompose the \sC dynamical system $(\mA,\rho)$ into the direct sum $\oplus_i (\mA_i,\rho_i)$, where $\mA_i := \chi_i \cdot \mA$ is the \sC algebra reduced by $\chi_i$, $\rho_i (\chi_i a) := \chi_i \rho (a)$, $a \in \mA$. Each endomorphism $\rho_i$ carries a weak permutation symmetry $\eps_i : p \mapsto \chi_i \eps (p)$ if $\rho$ does. Following Doplicher and Roberts \cite[\S 4]{DR87}, we may say that a permutation symmetry of $\rho$ has dimension $d \in \bN$ if the kernel of the morphism $C^*(\bP_\infty) \ra \mA$ is the ideal generated by the totally antisymmetric projection in $C^* (\bP_{d+1})$. But by the above considerations, a more general situation may arise, namely that every permutation symmetry $\eps_i$ has a distinct dimension $d_i$.

\begin{defn}
\label{dim_pm}
Let $\rho$ be a unital endomorphism of a \sC algebra $\mA$ carrying weak permutation symmetry, $d : \spzro \ra \bN$ a continuous (locally constant) map with range a finite set $\left\{ d_i \right\}$ of positive integers. We say that $\rho$ has rank $d$ if every endomorphism $\rho_i \in {\bf end}\mA_i$ has dimension $d_i$.
\end{defn}

Thus, the rank $d$ defines an element of the sheaf cohomology group $H^0(\spzro,\bZ)$ (i.e. a locally constant $\bZ$-valued map on $\spzro$). If a vector bundle $\mE \ra X$ has a nonconstant rank, then the canonical endomorphism of the associated Cuntz-Pimsner algebra has permutation symmetry carrying the same rank as $\mE$. In order to simplify the esposition, in the sequel we will restrict our attention to the case $d \in \bN$ (with $d > 1$). Note in fact that there is no loss of generality, thanks to the above direct sum decomposition of $(\mA,\rho)$.

\begin{lem}
\label{lem_da}
Let $\rho$ be an endomorphism of a unital \sC algebra $\mA$ with a rank $d$ weak permutation symmetry $p \mapsto \eps(p)$, $p \in \bP_\infty$. Then, for every rank $d$ vector bundle $\mE \ra \spzro$ there is a unique dual action $\mu : \wa \mUE \ra \wa \rho$, with $\mu (\theta) = \eps$, $\mu \circ \sigma_\mUE = \rho \circ \mu$. Moreover, suppose there is a Hilbert $\zro$-bimodule $\mR \subseteq ( \iota , \rho^d )_\eps$ with $\mR \mR^* = \zro P_{\eps,d}$; then, for every rank $d$ vector bundle $\mE \ra \spzro$ with a bimodule isomorphism $\beta : (\iota , \lambda \mE) \ra \mR$, the dual action $\mu$ uniquely extends to a dual $\mSUE$-action, by requiring $\mu (y) := \beta (y)$, $y \in (\iota,\lambda \mE)$.
\end{lem}

\begin{proof}
By \cite[Cor.4.4]{DR87}, there is a dual action $\wa {\ud} \ra \wa \rho$, which extends in the natural way to a dual action $\mu_0 : \wa{\mUE} \ra \wa \rho$, i.e. a $\zro$-morphism $\mu_0 : \mO_{\mUE} \ra \mA$ (recall $\mO_{\mUE} \simeq \zro \otimes \mO_{\ud}$). By Lemma \ref{cp_mu}, we know that there is a $\zro$-morphism $\mu : \mO_{\mSUE} \ra \mA$ extending $\mu_0$, and such that $\mu (\iota , \lambda \mE) = \mR$. In order to prove that $\mu$ is actually a dual action, we have to verify that 
\begin{equation}
\label{cp_int}
\mu \circ \sigma_{\mSUE} = \rho \circ \mu \ .
\end{equation}
\noindent Now, observe that since $\mu_0$ is a dual action, then (\ref{cp_int}) is verified for elements of $(\sigma^r_{\mSUE} , \sigma^r_{\mSUE}) = (\mE^r , \mE^r)_{\mSUE}$ (in fact, $(\mE^r , \mE^r)_{\mSUE} = (\mE^r , \mE^r)_{\mUE}$, see \cite[Cor.4.18]{Vas04}). Since $\mO_\mSUE$ is generated by $\mO_\mUE$ and $(\iota , \lambda \mE)$, it remains to verify (\ref{cp_int}) only for elements of $(\iota , \lambda \mE)$. By hypotesis $\mR \subset (\iota , \rho^d)_\eps$, so that if $R \in \mR$, then $\eps(d,1)R = \rho (R)$; now, if $y \in (\iota , \lambda \mE)$, then $\theta(d,1) y = \sigma_{\mSUE} (y)$ and $\mu (\theta(d,1) y) = \eps (d , 1) \mu (y) = \mu \circ \sigma_{\mSUE} (y)$. But $\mu (y) \in \mR \subset (\iota , \rho^d)_\eps$, so that $\eps(d,1) \mu (y) = \rho \circ \mu (y)$ and $\mu \circ \sigma_{\mSUE} (y) = \rho \circ \mu (y)$. We now prove that $\mu (\sigma^r_{\mSUE} , \sigma^{r + kd}_{\mSUE}) \subseteq (\rho^r , \rho^{r + kd})$. When $k = 0$, we know that such an inclusion is verified, $\mu_0$ being a dual action. Let $\left\{ S_L \right\} \subset ( \iota , \lambda \mE^k)$ be a subset with $\sum_L S^*_L S_L = 1$ (such a set is constructed by using a standard argument, involving the Serre-Swan Theorem and partitions of unity). Now, if $y \in (\sigma^r_{\mSUE} , \sigma^{r + kd}_{\mSUE})$, then $\mu (y) = \sum_L \mu_0(y \sigma^r_{\mSUE}(S^*_L)) \cdot \mu \circ \sigma^r_{\mSUE}(S_L)$, with $\mu_0(y \sigma^r_{\mSUE}(S^*_L)) \subseteq (\rho^{r+kd} , \rho^{r+kd})$, and $\mu \circ \sigma^r_{\mSUE}(S_L) = \rho^r \circ \mu (S_L) \subseteq \rho^r (\mR^k) \subset (\rho^r , \rho^{r+kd})$.
\end{proof}

We now introduce the notion of special conjugate property in our general setting. Note that in the case $\mZ = \bC$, such a notion is implemented by an isometry (see Def.\ref{def12}), which can be identified as a generator of the totally antisymmetric tensor power of the Hilbert space inducing $\rho$, if we construct the crossed product \cite[Thm.4.1]{DR89A}. In analogy with such a construction, we define the special conjugate property in terms of a Hilbert bimodule in $\mA$, corresponding to the module of continuous sections of the totally antisymmetric bundle.

\begin{defn}[Special Conjugate Property]
\label{def13}
Let $\rho$ be a unital endomorphism of a \sC algebra $\mA$ with weak permutation symmetry $\eps : \bP_\infty \ra \mA$. We say that $\rho$ satisfies the special conjugate property if for some $d \in \bN$, $d > 1$, there is a finitely generated Hilbert $\zro$-bimodule $\mR \subset ( \iota,\rho^d )_\eps$ (i.e. $\eps(d,1)R = \rho(R)$, $R \in \mR$) such that
\begin{equation}
\label{scp1}
R^* \rho (R') \ = \ (-1)^{d-1} d^{-1} R^* R' \ , \ R,R' \in \mR \ ;
\end{equation}
\begin{equation}
\label{scp2}
{\mR}{\mR}^*  \ := \ 
{\mathrm{closed \ span}}\left\{ R' R^*  :  R , R' \in \mR \right\} \ = \ 
\zro P_{\eps,d} \ .
\end{equation}
\end{defn}

The condition (\ref{scp2}) forces $\mR$ to be the module of sections of a line bundle over the spectrum $\spzro$. The condition (\ref{scp1}) implies that for every set $\left\{ R_i \right\}$ of normalized generators of $\mR$, the equation
\begin{equation}
\label{eq_gen}
\sum \limits_i R_i^* \rho (R_i) = (-1)^{d-1} d^{-1} 1
\end{equation}
\noindent holds, in analogy with (\ref{eq_scp1}). Note that we could more generally state our special conjugate property Def.\ref{def13} by considering a non-constant rank $d \in H^0(\spzro,\bZ)$, by performing the decomposition $(\mA,\rho) \simeq \oplus_i (\mA_i , \rho_i)$ according to Def.\ref{dim_pm}, and by giving the previous definition for each index $i$. The result is a Hilbert $\zro$-bimodule $\mR$ in $\mA$ with a natural decomposition $\mR \simeq \oplus_i \mR_i \subset \oplus_i (\iota , \rho_i^{d_i})_\eps$, where each $\mR_i$ satisfies (\ref{scp1}, \ref{scp2}) for a fixed dimension $d_i$.

\begin{lem}
\label{lem_uni_r}
There is only one $\zro$-bimodule $\mR \subset (\iota,\rho^d)$ satisfying (\ref{scp1},\ref{scp2}).
\end{lem}

\begin{proof}
Let $\mathcal S \subset (\iota,\rho^d)$ be another $\zro$-module satisfying Def.\ref{def13}, with a set of generators $\left\{ S_j \right\}$ (i.e., $\sum_j S_j S_j^* = P_{\eps,d}$). We consider a set of normalized generators $\left\{ R_i \right\}$ for $\mR$ (i.e., $\sum_i R_i R_i^* = P_{\eps,d}$, $\sum_i R_i^* R_i = 1$), and define $z_{ij} := R_i^* S_j \in \mZ$ for each $i,j$. In order to prove the lemma, it will suffice to show that $z_{ij} \in \zro$. For this purpose, observe that for each $R \in \mR$ the identity $R_i^* \rho (R) = (-1)^{d-1} d^{-1} R_i^* R$ holds; so that, by multiplying on the left by $R_i$, and by summing over $i$, we find 
\begin{equation}
\label{exp_rho}
P_{\eps,d} \rho (R) = (-1)^{d-1} d^{-1} R \ .
\end{equation}
\noindent In the same way, for each $S \in \mathcal S$, we obtain $P_{\eps,d} \rho (S) = (-1)^{d-1} d^{-1} S$. Now, observe that $R_k R_i^* = f_{ik} P_{\eps,d}$, where $f_{ik} \in \zro$, so that $R_k z_{ij} = z_{ij} R_k = R_i^* S_j R_k = f_{ik} P_{\eps,d} S_j = f_{ik} S_j$. Thus, $\rho (R_k) \rho(z_{ij}) = f_{ik} \rho (S_j)$; by multiplying this last identity on the left by $P_{\eps,d}$, and by applying (\ref{exp_rho}) for $R_k,S_j$, we obtain
\[
(-1)^{d-1} d^{-1} R_k \rho(z_{ij}) = (-1)^{d-1} d^{-1} f_{ik} S_j 
= (-1)^{d-1} d^{-1} R_k z_{ij} \ ,
\]
\noindent and $\rho(z_{ij})^* R_k^* = z_{ij}^* R_k^*$. By multiplying on the right by $R_k$, and by summing over $k$ we obtain $\rho(z_{ij})^* = z_{ij}^*$.
\end{proof}

The first Chern class in $H^2(\spzro , \bZ)$ (sheaf cohomology) of the line bundle having $\mR$ as module of continuous sections is a complete invariant of $\mR$; the unicity established in Lemma \ref{lem_uni_r} implies that it is actually an invariant of $\rho$. We denote it by $c_1 (\rho)$.

We now give an application of \cite[Thm.4.3]{DR87}, and compute the rank of a weak permutation symmetry $p \mapsto \eps(p)$ for an endomorphism $\rho$ satisfying the special conjugate property Def.\ref{def13}. \cite[Thm. 4.3]{DR87} asserts that if $\phi : \mA \ra \mA$ is a left inverse for $\rho$ with $\phi (\eps) = \lambda 1$, $\lambda \in \bC$, then $d$ is recovered by the formula $\lambda = \pm d^{-1}1$, and $\lambda \in \left\{ 0 \right\} \cup \left\{ \pm k^{-1} , k \in \bZ - \left\{ 0 \right\} \right\}$.

\begin{lem}
\label{dim_perm}
Let $\rho$ be an endomorphism of a unital \sC algebra $\mA$ satisfying the special conjugate property Def.\ref{def13}. Then $\rho$ has weak permutation symmetry of rank $d$.
\end{lem}

\begin{proof}
We consider a normalized finite set $\left\{ R_i \right\}$ of generators of $\mR$, and define the left inverse
\[
\phi (a) := \sum \limits_i R_i^* \rho^{d-1} (a) R_i  , \quad a \in \mA \ .
\]
\noindent In order to apply \cite[Thm.4.3]{DR87}, we have to explicitly compute $\phi (\eps)$. For this purpose, note that $\bS^{d-1} (1,1) = (d-1,1)(d-1)$, so that $\rho^{d-1}(\eps) = \eps(d-1,1) \eps(d,1)$. Thus, $\phi(\eps) = \sum_i R_i^* \eps(d-1,1) \eps(d,1) R_i$. Now, $\eps(d,1) R_i = \rho (R_i)$ (in fact, $\mR \subset (\iota,\rho^d)_\eps$), and $\eps(1,d-1) P_{\eps,d} = \mathrm{sign} (1,d-1) P_{\eps,d} = (-1)^{d-1} P_{\eps,d}$, so that $\eps(1,d-1) R_i = (-1)^{d-1} R_i$. Thus, $\phi(\eps) = (-1)^{d-1} \sum_i R_i^* \rho (R_i)$. By applying (\ref{eq_gen}), we find $\phi(\eps) = d^{-1} 1$.
\end{proof}

We recall the reader to the notation (\ref{def_hzro}).

\begin{defn}
\label{cro}
A unital endomorphism $\rho$ of a \sC algebra $\mA$ is {\bf weakly special} if satisfies the special conjugate property Def.\ref{def13}. The class of $\rho$ is
\[
c(\rho) := c_0(\rho) \oplus c_1(\rho) \in H^{0,2}( \spzro , \bZ ) \ ,
\]
\noindent where $c_0(\rho)$ is the rank of the weak permutation symmetry, and $c_1(\rho)$ is the Chern class of the Hilbert $\zro$-bimodule $\mR$ introduced in Def.\ref{def13}. Moreover, $\rho$ is said {\bf quasi-special} if the permutation quasi-symmetry Def.\ref{def_gps} holds.
\end{defn}

Note that the classical special conjugate property Def.\ref{def12} is satisfied if and only if $c_1(\rho)$ vanishes, and $c_0(\rho)$ is a constant map. In that case, $\mR$ is generated as a Hilbert $\zro$-bimodule by a partial isometry with support $P_{\eps , d}$.

\begin{ex} {\em 
Let $\mE \ra X$ be a rank $d$ vector bundle, $G \subseteq \mSUE$. Then, $\sigma_G$ is quasi-special in the sense of Def.\ref{cro}. In fact, $\sigma_G$ has permutation symmetry and satisfies the special conjugate property Def.\ref{def13}, implemented by the module of continuous sections of $\lambda \mE$ (see Sec.\ref{apdx}). Thus, $c(\sigma_G) = d \oplus c_1(\mE)$, where $c_1(\mE)$ is the first Chern class of $\mE$. }
\end{ex}

\begin{ex} {\em 
We retain the notation of Ex.\ref{ps_bp}. The shift endomorphism $\tau_G \in {\bf end}_X \cmg$ has weak permutation symmetry; furthermore, if $G \subseteq \mSUE$, we find that $\tau_G$ satisfies the special conjugate property with $\mR = (\iota , \lambda \mE) \subseteq (\iota , \tau_G^d)$. So that $\tau_G$ is a weakly special endomorphism in the sense of Def.\ref{cro}. }
\end{ex}

\begin{ex} {\em 
Let $(\mF , G)$ be a minimal regular Hilbert \sC system in the sense of \cite{BL04}, with fixed-point algebra $\mA$. Canonical endomorphisms of $\mA$ in the sense of \cite[Def.3.6]{BL04} are quasi-symmetric (\cite[Thm.4.10]{BL04}). Moreover, every canonical endomorphism is dominated by a quasi-special endomorphism $\rho$ (see \cite[Thm.3.4]{DR89}), satisfying Def.\ref{def13} with $\mR = \zro \cdot R$ for a suitable partial isometry $R \in (\iota , \rho^d)_\eps$. $R$ is a determinant of $\rho$ in the sense of \cite[\S 3]{DR89}. Thus, $c_1 (\rho) = 0$. }
\end{ex}

\subsection{Weakly special endomorphisms, and duality.}

Let $\mZ$ be a \sC algebra, $\mM$ a Hilbert $\mZ$-bimodule, $G$ a closed group of unitary {\em $\mZ$-bimodule} operators of $\mM$ (i.e., the elements of $G$ commute with the left $\mZ$-module action). Then, every internal tensor power $\mM^r$, $r \in \bN$, is a $G$-Hilbert bimodule in the natural way $g , \psi \mapsto g_r \psi$, $\psi \in \mM^r$, $g \in G$, $g_r := g \otimes \cdots \otimes g$. We denote by ${\wa G}_\mM$ the \sC category with objects $\mM^r$, $r \in \bN$, and arrows the spaces $(\mrs)_G$ of $G$-equivariant right $\mZ$-module operators from $\mM^r$ into $\mM^s$. Let now $X$ be a compact Hausdorff space, $\mZ$ a $C(X)$-algebra, and $\mM \simeq \wE \otimes_{C(X)} \mZ$ a nc-pullback of a vector bundle $\mE \ra X$. If $G \subseteq \mUE$ is a closed group, then $\mM$ is a $G$-Hilbert $\mZ$-bimodule, and there is an obvious inclusion $\wa G \subseteq {\wa G}_\mM$ (i.e., $(\ers)_G \subseteq (\mrs)_G$, $r,s \in \bN$). For details about this construction, see \cite[\S 5]{Vas04}.

\begin{prop}
\label{cor_dual_1}
Let $\mE \ra X$ be a vector bundle, $G \subseteq \mSUE$, $\mc$ the crossed product by a dual action $\mu : \wa G \ra \wa \rho \subset \tend$; then, $\rho$ is weakly special. Suppose $\mA$ is the fixed-point algebra w.r.t. the $G$-action on $\mc$, and $\mA' \cap (\mc) = \mA \cap \mA' =: \mZ$; then, $\mu$ induces isomorphisms of tensor \sC categories
\begin{equation}
\label{iso_rhoeps}
\wa G \stackrel{\simeq}{\lra} \wa \rho_\eps \ \ : \ \ 
(\ers)_G \simeq (\rhors)_\eps \ \ ,
\end{equation}
\begin{equation}
\label{iso_rho}
{\wa G}_\mM \stackrel{\simeq}{\lra} \wa \rho \ \ : \ \ 
(\mrs)_G \simeq (\rhors) \ \ ,
\end{equation}
\noindent $r,s \in \bN$, where $\mM$ is the $G$-Hilbert $\mZ$-bimodule defined by (\ref{def_bm}).
\end{prop}

\begin{proof}
Recall that $\mM = \wE \cdot \mZ$ (Prop.\ref{cp_mcr}). For every $r,s \in \bN$, there are natural identifications $(\mrs) = {\mathrm{span}} \left\{ \psi_L z \psi_M^*  , z \in \mZ  \right\}$, where $\psi_L \in (\iota , \mE^s)$, $\psi_M \in (\iota , \mE^r)$ are defined by (\ref{def_psiL}). Since $(\iota , \mE^r) \subseteq (\iota , \sigma_\mE^r)$ (see Rem.\ref{rem_ce}), we find 
\begin{equation}
\label{eq_rs}
\psi_L z \psi_M^* \rho^r(a) = \psi_L z a \psi_M = \rho^s (a) \psi_L z \psi_M^*
\ \ , \ \ 
a \in \mA \ .
\end{equation}
\noindent If $t \in (\mrs)_G$, then by hypothesis $t \in \mA$; moreover, (\ref{eq_rs}) implies that $t \in (\rhors)$. Viceversa, if $t \in (\rho^r,\rho^s)$, then $t_{LM} := \psi_L^* t \psi_M \in \mA' \cap (\mc) = \mZ$, and $t = \sum_{L,M} \psi_L t_{LM} \psi_M^* \in (\mrs)$. This proves that $(\mrs)_G = (\rhors)$, $r,s \in \bN$, thus we proved (\ref{iso_rho}). 

Now, Rem.\ref{thm_bgae} implies that $\rho$ has weak permutation symmetry. Since $G \subseteq \mSUE$, then $(\iota , \lambda \mE) \subseteq (\iota , \mE^d)_G$, and $\mR := \mu (\iota , \lambda \mE) \subseteq (\iota , \rho^d)_\eps$; thus, $\rho$ is weakly special. By Rem.\ref{rem_si_da}, we find $\mu (\ers)_G \subseteq (\rhors)_\eps$, $r,s \in \bN$. Viceversa, if $t \in ( \rhors )_\eps$, then with the above notation $\rho (t_{LM}) = \sigma_\mE (\psi_L^*) \rho (t) \sigma_\mE (\psi_M) = \sigma_\mE (\psi_L^*) \eps (s,1) t \eps (1,r) \sigma_\mE (\psi_M) = t_{LM}$; in fact, 
\begin{equation}
\label{eq_spsi}
\sigma_\mE (\psi_L) = \mu \circ \theta(s,1) \psi_L = \eps (s,1) \psi_L
\end{equation}
\noindent (see \cite[(4.16)]{Vas04}). Thus, $t_{LM} \in \zro$, and this implies $t \in (\ers)$. Since $t$ is $G$-invariant, we proved (\ref{iso_rhoeps}). 

\end{proof}

The condition that $\mA$ is the fixed-point algebra w.r.t. the $G$-action is automatically verified if $G$ is compact (see remarks after Thm.\ref{thm38}), or locally trivial in the sense of next subsection (see Lemma \ref{lem38}).

\subsection{Locally trivial groups, and a pullback construction.}
\label{sec_pb}

In the sequel we will use the notation $H_d := \bC^d$, so that every $(H_d^r , H_d^s)$ is identified with the matrix space $\bM_{d^r , d^s}$, $r,s \in \bN$. Note that if $G_0 \subseteq \ud$, then the canonical action (\ref{eq_def_G^x}) is defined; as usual, we denote by $(H_d^r , H_d^s )_{G_0}$ the corresponding $G_0$-invariant vector spaces.

\

A closed group $G \subseteq \mUE$ is said {\em locally trivial} if the bundle $\mG \ra X$ considered in Sec.\ref{apdx} is locally trivial (see \cite[Def.4.11]{Vas04}). In such a case, we find that every $G^x$, $x \in X$, coincides with a fixed compact Lie group $G_0 \subseteq \ud$ (for example, $G := \mSUE$ is locally trivial with $G^x \equiv \sud$). If $\left\{ U_i \subseteq X \right\}_{i}$ is a finite open cover trivializing $\mG$, then for every index $i$ there is an isomorphism $\pi_i : \mG |_{U_i} \ra U_i \times G_0$; since by definition $G$ is the group of continuous sections of $\mG$, every $g \in G$ may be described by a family of continuous maps $g_i : U_i \ra G_0$, with $( x , g_i(x)) := \pi_i \circ g(x)$, $x \in U_i$. A locally trivial group is typically noncompact (for example, if $\mE := X \times \bC^d$, then $G := \mUE$ is the 'loop group' of continuous maps from $X$ into $\ud$).

\begin{lem}
\label{lem_gcut}
Let $G_0$ be connected. For every $x_0 \in X$, there is a closed neighborhood $W$ of $x_0$, and an open $U \subseteq X$, such that: 
\begin{enumerate}
\item  $W$ is strictly contained in $U$, and $U$ trivializes $\mG$ (i.e. $\mG |_U \simeq U \times G_0$);
\item  for every $u \in G_0$, there is a non-unique $\wt u \in G$ such that $\wt u \ |_W \equiv u$, $\wt u \ |_{X- U } \equiv 1$ (in particular, $\wt u (x_0) = u$).
\end{enumerate}
\end{lem}

\begin{proof}
Let $x_0 \in X$, and $U \subseteq X$ be an open neighborhood of $x_0$ trivializing $\mG$, so that $\mG |_U \simeq U \times G_0$. We consider a closed neighborhood $W$ of $x_0$ strictly contained in $U$. Moreover, we consider a cutoff $\lambda : X \ra [0,1]$ such that $\lambda |_W \equiv 1$, $\lambda_{X - U} \equiv 0$. Let $u \in G_0$, and $\gamma : [0,1] \ra G_0$ a path such that $\gamma (0) = 1$, $\gamma (1) = u$. We define $\wt u \in G$, $\wt u (x) := \gamma \circ \lambda (x)$, $x \in X$. Note that $\wt u \ |_W \equiv u$, $\wt u \ |_{X-U} \equiv 1$.
\end{proof}

For the rest of the present subsection, we will consider a locally trivial group $G \subseteq \mSUE$ with fibre $G_0 \subseteq \sud$, a unital \sC algebra $\mA$, and a fixed dual action $\mu : \wa G \ra \wa \rho$, $\rho \in \tend$. In order to simplify the notation, we write $X := \spzro$.

\begin{lem}
\label{lem38}
Suppose that $G_0$ is connected. Then, the fixed-point algebra w.r.t. the action $\alpha : G \ra {\bf aut}_X (\mc)$ is $\mA \otimes_\mu 1 \simeq \mA$.
\end{lem}

\begin{proof}
It is clear that $\mA \subseteq (\mc)^\alpha$. In order to prove the opposite inclusion, for every $x \in X$ we consider a closed neighborhood $W_x$ as in Lemma \ref{lem_gcut}, in such a way that $\left\{ \right. \dot{W}_x \subseteq X \left. \right\}_x$ is an open cover ($\dot{W}_x$ denotes the interior). By compactness, there is a finite subcover $\left\{ \right. \dot{W}_k \left. \right\}_{k=1}^n$; we consider a subordinate partition of unity $\left\{ \lambda_k \right\}_{k=1}^n$. We denote by $\mA_k := C_0(\dot{W_k}) \mA :=$ closed span$\left\{ \right. f a , f \in C_0(\dot{W_k}) , a \in \mA \left. \right\}$ the restriction of $\mA$ over $\dot{W}_k$ as an upper semicontinuous bundle; in the same way, we define $\mB_k := C_0(\dot{W_k}) (\mc)$. It is clear that $\mA_k \subset \mB_k$ are $C_0(\dot{W_k})$-algebras. For every $k = 1, \ldots , n$, we define the following action:
\[
\alpha_k : G_0 \ra {\bf aut}_{ \dot{W_k} } \mB_k 
\ \ , \ \ 
\alpha_k (u) := \alpha_{\wt u} \ \ ,
\]
\noindent where $\wt u \in G$ is defined as in Lemma \ref{lem_gcut}, and $\alpha_{\wt u}$ is defined according to the action introduced in Thm.\ref{thm38}. Since $\wt u |_{W_k} \equiv u$, the above action does not depends on the particular choice of $\wt u$. Now, the following $\mA_k$-valued invariant mean is defined on $\mB_k$: $$m_k (at) := a \int_{G_0} \alpha_{\wt u} (t) du = a \int_{G_0} \wa u (t) du \ ,$$ $a \in \mA_k$, $t \in C_0(\dot{W}_k) (\ers)$ (in fact, by Lemma \ref{def_ra}, the Haar mean of $t$ belongs to $C_0(\dot{W}_k) \pi_{W_k} (\ers)_G \subset \mA_k$). By using such invariant mean, it is easy to verify that the fixed-point algebra w.r.t. the $\alpha_k$-action is given by $\mA_k$, as in \cite[Lemma 2.8, Thm.3.2]{DR89A}. If $b \in (\mc)^\alpha$, then $\lambda_k b \in \mB_k$, and $[\alpha_k (u)] (\lambda_k b) = \alpha_{\wt u} (\lambda_k b)  = \lambda_k b$; thus, $\lambda_k b \in \mA_k$. Since $b = \sum_k \lambda_k b_k$, we conclude that $b \in \mA$.
\end{proof}

It is natural to conjecture that the fixed-point algebra coincides with $\mA$ also in the general case (i.e., $G$ not locally trivial). Since our main application will concerns the case $G = \mSUE$, for our purposes it will suffice to consider locally trivial groups.

\

We denote by $\mC := C(\Omega)$ the centre of $\mc$; it is clear that $\mC$ is a $C(X)$-algebra. Now, the following properties hold:

\begin{enumerate}

\item  There is a surjective map $p : \Omega \ra X$ (i.e., the adjoint map of the inclusion $C(X) \hra \mC$). Thus, $\mC$ is an upper-semicontinuous bundle $(\mC , \eta_x : \mC \ra \mC_x)$, with fibres the abelian \sC algebras $\mC_x \simeq C(p^{-1}(x))$, $x \in X$. The image of $c \in \mC$ w.r.t. $\eta_x$ is simply the restriction of $c$ over $p^{-1}(x)$.

\item  Let $g(x) \in G_0$ denote the evaluation of $g \in G$ over $x \in X$. $G$ acts by $C(X)$-automorphisms on $\mC$, in such a way that the fixed-point algebra is $C(X)$; so that, there are automorphic actions $\alpha_x : G_0 \ra {\bf aut}\mC_x$, in such a way that $\eta_x \circ \alpha_g  = \alpha_{ x , g(x) } \circ \eta_x$, $g \in G$ (see (\ref{eq_intact})).

\end{enumerate}

\begin{lem}
\label{lem_fpa}
Suppose that $G_0$ is connected. Then, the fixed-point algebras w.r.t. the actions $\alpha_x : G_0 \ra {\bf aut} \mC_x$, $x \in X$, are $\mC_x^{\alpha_x} \simeq \bC$.
\end{lem}

\begin{proof}
Let $v \in X$, $c_0 \in \mC_v^{\alpha_v}$; since $p^{-1}(v)$ is closed in $\Omega$, by the Tietze theorem there is $c_1 \in \mC$ such that $\eta_v (c_1) = c_0$. We now consider closed neighborhoods $W \subset U$, $W \ni v$, as in Lemma \ref{lem_gcut}; note that we could assume that $c_1$ has support contained in $p^{-1}(W)$ (otherwise, we multiply by a cutoff $\lambda : X \ra [0,1]$, $\lambda (v) = 1$, $\lambda |_{X - W} \equiv 0$). Let $\mC_W := C_0(\dot{W}) \mC \subset \mC$ denote the restriction of $\mC$ over the interior $\dot{W}$; we define the automorphic action
\[
\alpha_W : G_0 \ra {\bf aut} \mC_W 
\ \ , \ \ 
[\alpha_W (u)] (c) := \alpha_{\wt u} (c) \  ,
\]
\noindent where $\wt u \in G$ is assigned to every $u \in G_0$ as in Lemma \ref{lem_gcut}. Since $\wt u (x) \equiv u$, $x \in W$, it turns out that $\alpha_W (u)$ does not depends on the choice of $\wt u$. It is clear that $\eta_x \circ \alpha_W (u) = \eta_x \circ \alpha_g = \alpha_{x,u} \circ \eta_x$, $x \in W$, $g \in G$, $u := g(x) \in G_0$. We define $$c_2 := \int_{G_0} \alpha_{\wt u} (c_1) \ du  \ \in \ \mC_W \subset \mC \ \ .$$ By definition, $c_2$ is $\alpha_W$-invariant, i.e. $\alpha_{\wt u} (c_2) = c_2$, $u \in G_0$, and $\eta_x (c_2) = c_0$; anyway, more in general, we find
\[
[\alpha_g (c_2)] (\omega) = 
[\eta_x \circ \alpha_g (c_2)] (\omega) =
[\alpha_{\wt u} (c_2) ] (\omega) =
c_2(\omega) \ \ ,
\]
\noindent $\omega \in p^{-1}(W)$, $x := p(\omega)$, $g \in G$, $u := g(x)$. Since $c_2$ has support contained in $p^{-1}(W)$, we conclude that $c_2$ is $\alpha$-invariant, i.e. $c_2 \in C(X)$. But this means that $c_2$ is constant on each fibre $p^{-1}(x)$, $x \in X$; in particular, $\eta_v (c_2) = c_0$ is constant on $p^{-1}(v)$.
\end{proof}

Thus, $p^{-1}(x)$ is (noncanonically) homeomorphic to the homogeneous space $\left. L_x \right/ G_0$, where $L_x \subseteq G_0$ is the stabilizer of some (noncanonical) $\omega_x \in p^{-1}(x)$ (\cite{EH}). By functoriality, every $\alpha_x$ induces an ergodic $G$-action 
\[
p^{-1} (x) \times G \ra p^{-1} (x)
\ \ , \ \ 
\omega , g \mapsto \omega \cdot g := \omega \circ \alpha_{ x , g(x) } \ .
\]
\noindent Note that $\omega \cdot g$ actually depends only on $g(x) \in G_0$. For every $\omega \in p^{-1}(x)$, the orbit $\left\{ \omega \cdot g , g \in G \right\}$ coincides with $p^{-1}(x)$. We can now perform the following construction.

\begin{enumerate}

\item  We denote by $p_* \mE := \mE \times_X \Omega \ra \Omega$ the pullback of $\mE \ra X$ over $\Omega$. The $C(\Omega)$-bimodule $p_* \wE$ of continuous sections of $p_* \mE$ is isomorphic to $\wE \otimes_{C(X)} C(\Omega)$; more in general, $(p_* \mE^r , p_* \mE^s) \simeq (\ers) \otimes_{C(X)} C(\Omega)$. Let $\pi_{rs} : \mE^{r,s} \ra X$ be the vector bundle having $(\ers)$ as module of continuous sections. Then, $(p_* \mE^r , p_* \mE^s)$ is the module of continuous sections of $p_* \mE^{r,s} := \mE^{r,s} \times_X \Omega \ra \Omega$. Every $t \in (p_* \mE^r , p_* \mE^s)$ can be regarded as a continuous map $t : \Omega \ra \mE^{r,s}$ such that $p = \pi_{rs} \circ t$.

\item  Every $p_* \mE^{r,s}$ is a $G$-vector bundle in the sense of \cite[\S 1.6]{Ati}, with {\em $G$ acting nontrivially on $\Omega$}. We have the $G$-actions
\begin{equation}
\label{def_gact}
\left\{
\begin{array}{ll}
G \times  (p_* \mE^r , p_* \mE^s) \ra (p_* \mE^r , p_* \mE^s)    \\
g,t \mapsto g \times t \ : \ 
(g \times t) (\omega) := 
g_s (x) \cdot t (\omega \cdot g^*) \cdot g_r^* (x) \ \ ,
\end{array}
\right.
\end{equation}
\noindent $x := p (\omega)$, $g_r := g \otimes \cdots \otimes g \in (\mE^r , \mE^r)$. We define $p_* \wa G$ as the tensor \sC category with objects the tensor powers of $p_* \mE$, and  sets of arrows the $C(X)$-bimodules $(p_* \mE^r , p_* \mE^s)_G$ of $G$-invariant elements in $(p_* \mE^r , p_* \mE^s)$. Let $t \in (p_* \mE^r , p_* \mE^s)$; then,  $t \in (p_* \mE^r , p_* \mE^s)_G$ if and only if
\begin{equation}
\label{rel_inv}
t (\omega \cdot g) =
g_s (x) \cdot t (\omega) \cdot g^*_r (x) \ \ ,
\end{equation}
\noindent $g \in G$, $\omega \in \Omega$, $x := p(\omega)$; i.e., $t = g \times t$, $g \in G$.

\end{enumerate}

Let now $W \subseteq X$ be a closed set; then, $p^{-1} (W) \subseteq \Omega$ is closed and $G$-stable. We consider the restriction $p_* \mE_W := p_* \mE |_{p^{-1}(W)}$; in the same way, we define $p_* \mE^{r,s}_W := p_* \mE^{r,s} |_{p^{-1} (W)}$. For every $r,s \in \bN$, there is a natural epimorphism
\[
\pi_W : ( p_* \mE^r , p_* \mE^s) \ra ( p_* \mE^r_W ,  p_* \mE^s_W ) \ ,
\]
\noindent defined as the restriction of $( p_* \mE^r , p_* \mE^s)$ over $p^{-1}(W)$ as a continuous field of Banach spaces (\cite[10.1.7]{Dix}). If $p^{-1}(W)$ trivializes $p_* \mE$, then 
\[
\left\{
\begin{array}{ll}
p_* \mE^{r,s}_W \simeq p^{-1} (W) \times ( H_d^r , H_d^s ) 
\\
( p_* \mE_W^r , p_* \mE_W^s) = C( p^{-1} (W) ) \otimes ( H_d^r , H_d^s ) 
\end{array}
\right.
\]

The following lemma expresses in rigorous terms the intuitive idea that the $G$-action on $p_* \mE$ locally looks like a $G_0$-action.

\begin{lem}
\label{def_ra}
Let $G_0$ be connected. Then, for every $x_0 \in X$ there is a closed neighborhood $W \ni x_0$ such that $p_* \mE_W$ is a $G_0$-vector vector bundle. Moreover,
\begin{equation}
\label{eq_def_ra}
( p_* \mE_W^r , p_* \mE_W^s  )_{G_0} = \pi_W (  p_* \mE^r , p_* \mE^s  )_G
\ \ , \ \ 
r,s \in \bN \ \ .
\end{equation}
\end{lem}

\begin{proof}
For every $u \in G_0$, we consider $\wt u$ as in Lemma \ref{lem_gcut}. Let $t \in (  p_* \mE^r , p_* \mE^s  )$; for every $\omega \in p^{-1} (W)$, by (\ref{def_gact}) we find
\[
(\wt u \times t) (\omega) =  
\wt u_s (x) \cdot t (\omega \cdot \wt u^*) \cdot \wt u^*_r (x) =
u_s \cdot t (\omega \cdot g^*) \cdot u_r^* \ .
\]
\noindent Thus, the vector field $\pi_W (\wt u \times t) \in ( p_* \mE_W^r , p_* \mE_W^s )$ does not depend on the choice of $\wt u$; we define the $G_0$-action
\begin{equation}
\label{def_gw}
\left\{
\begin{array}{ll}
G_0 \times (  p_* \mE_W^r , p_* \mE_W^s ) \ra (  p_* \mE_W^r , p_* \mE_W^s )
\\ 
u , \pi_W(t) \mapsto u \times \pi_W (t) := \pi_W (\wt u \times t)  \  .
\end{array}
\right.
\end{equation}
\noindent In particular, for $r=0$, $s=1$ the Serre-Swan theorem implies that $p_* \mE_W$ is a $G_0$-vector bundle.

Let now $t \in (  p_* \mE^r , p_* \mE^s  )_G$. Then, $\pi_W(\wt u \times t) = \pi_W(t)$ for every $u \in G_0$; since $\wt u (x) = u$, $x \in W$, we find $\pi_W (t) \in ( p_* \mE_W^r , p_* \mE_W^s )_{G_0}$. Viceversa, let $t_0 \in ( p_* \mE_W^r , p_* \mE_W^s )_{G_0}$. Let $V \subseteq X$ be a closed set with interior $\dot{V}$, and such that there are strict inclusions $W \subset \dot{V}$, $V \subset U$ (where $U$ is chosen as by Lemma \ref{lem_gcut}). By the previous argument, the analogue for $V$ of the action (\ref{def_gw}) is defined. By the Tietze theorem \cite[Lemma 1.6.3]{Ati}, there is $t_1 \in ( p_* \mE^r_V , p_* \mE^s_V )_{G_0}$ such that $\pi_W (t_1) = t_0$. We now consider a cutoff $\lambda : X \ra [0,1]$, $\lambda |_W \equiv 1$, $\lambda_{X - V} \equiv 0$, and define $t_2 := \lambda t_1$. It is clear that $t_2 \in ( p_* \mE^r , p_* \mE^s )$, and $\pi_W (t_2) = t_0$. Let us verify that $t_2 \in (  p_* \mE^r , p_* \mE^s  )_G$; for every $g \in G$, $\omega \in V$, $x := p(\omega)$, we find
\[
(g \times t_2) (\omega) = 
( \wt {g(x)} \times t_2 ) (\omega) =
\lambda ( \wt {g(x)} \times t_1 ) (\omega) =
\lambda t_1 (\omega) =
t_2 (\omega)
\ ,
\]
\noindent where $\wt {g(x)} \in G$ is constructed as by Lemma \ref{lem_gcut} (note $\wt {g(x)} |_V \equiv g(x) \in G_0$), and $\wt {g(x)} \times t_1$ is defined as in (\ref{def_gw}). Thus, $t_2 \in$ $(  p_* \mE^r ,$ $ p_* \mE^s  )_G$, and the lemma is proven.
\end{proof}

\begin{prop}
\label{rhoine}
Suppose that $G_0 \subseteq \sud$ is connected. Then, for every $r , s \in \bN$,
\begin{enumerate}
\item  there is a natural isomorphism $(p_* \mE^r , p_* \mE^s) \simeq (\sigma_\mE^r , \sigma_\mE^s)$; 
\item  there is a natural isomorphism ${(\rhors)_\eps} \simeq (p_* \mE^r , p_* \mE^s)_G$, i.e. there is an isomorphism of tensor \sC categories $\wa \rho_\eps \simeq p_* \wa G$.
\end{enumerate}
\end{prop}

\begin{proof}
\
\begin{enumerate}
\item  By definition (Rem.\ref{rem_ce}), it is clear that $(\iota , \sigma_\mE)$ is generated as a $\mC$-module by elements of $\wE$, so that $(\iota , \sigma_\mE) = \wE \cdot \mC$, and we obtain isomorphisms $(\iota , \sigma_\mE) = \wE \cdot \mC \simeq \wE \otimes_{C(X)} \mC \simeq (\iota , p_* \mE)$. Thus, the isomorphisms $(p_* \mE^r , p_* \mE^s) \simeq (\sigma_\mE^r , \sigma_\mE^s)$ trivially follow.
\item  Let $t \in (p_* \mE^r , p_* \mE^s)_G$; then $t = \sum_{LM} \psi_L c_{LM} \psi_M^*$, where $\psi_L$, $\psi_M$ are defined according to (\ref{def_psiL}), and $c_{LM} \in \mC$. Thus, $(p_* \mE^r , p_* \mE^s)_G \subseteq (\sigma_\mE^r , \sigma_\mE^s)$. By construction, the action (\ref{def_gact}) coincides with the $G$-action on $\mc$ for elements of $(p_* \mE^r , p_* \mE^s)$; by Lemma \ref{lem38}, we conclude $(p_* \mE^r , p_* \mE^s)_G$ $\subseteq \mA$, and $(p_* \mE^r , p_* \mE^s)_G \subseteq (\rhors)$. Moreover, by applying (\ref{eq_spsi}), we find $$\rho (t) = \sum_{LM} \eps (s,1) \psi_L c_{LM} \psi_M^* \eps (1,r) = \eps (s,1) t \eps (1,r) \ ,$$ thus $t \in (\rhors)_\eps$. Viceversa, if $t \in {(\rhors)_\eps}$ then $t = \sum_{L,M} \psi_L t_{LM} \psi_M^*$, where $t_{LM} := \psi_L^* t \psi_M$. We now prove that $t_{LM}$ belongs to $\mC$: first note that 
\[
t_{LM} a = 
\psi_L^* t \psi_M a = 
\psi_L^* t \rho^r(a) \psi_M = 
\psi_L^* \rho^s(a) t \psi_M=
a t_{LM} \ ,
\]
\noindent $a \in \mA$, thus $t_{LM} \in \mA' \cap (\mc)$. Furthermore, $$\sigma_\mE (t_{LM}) = \sigma_\mE (\psi_L^*) \rho (t) \sigma_\mE (\psi_M) = \sigma_\mE (\psi_L^*) \eps (s,1) t \eps (1,r) \sigma_\mE (\psi_M) = t_{LM} .$$ Since $\wE \subseteq (\iota , \sigma_\mE)$, the previous computation implies $\psi t_{LM} = \sigma_\mE (t_{LM}) \psi = t_{LM} \psi$, $\psi \in \wE$, i.e. $t_{LM}$ commutes with elements of $\wE$. Since $\mc$ is generated by $\mA$, $\wE$, we conclude that $t_{LM}$ belongs to $\mC$; thus, we proved that $t \in (\sigma_\mE^r , \sigma_\mE^s) \simeq (p_* \mE^r , p_* \mE^s)$. Since by definition $t$ is $G$-invariant, the proposition is proven.
\end{enumerate}
\end{proof}

\begin{thm}
\label{soro_field}
Let $\mu : \wa G \ra \wa \rho$ be a dual action, where $G \subseteq \mSUE$ is locally trivial with connected fibre $G_0 \subseteq \sud$. Then, for every $x \in X$ there is a compact Lie group $L_x \subseteq G_0$, unique up to conjugation in $G_0$, such that:

\begin{enumerate}

\item  every ${(\rhors)_\eps}$ is a continuous field of Banach spaces, having fibres isomorphic to $( H_d^r , H_d^s )_{L_x}$, $x \in X$;

\item $\soro$ is a continuous bundle of \sC algebras, with fibres isomorphic to $\mO_{L_x}$, $x \in X$.

\end{enumerate}
\end{thm}

\begin{proof}
By the previous proposition, we identify ${(\rhors)_\eps}$ with $(p_* \mE^r , p_* \mE^s)_G$. Thus, every $t \in {(\rhors)_\eps}$ may be regarded as a continuous map $t : \Omega \ra \mE^{r,s}$ satisfying (\ref{rel_inv}); in particular, every restriction $t |_{p^{-1}(x)}$ takes values in the fibre $\mE^{r,s}_x \simeq ( H_d^r , H_d^s )$. For every $x \in X$, let us consider the vector space 
\[
\rho^{r,s}_x := 
\left\{ 
t_x := t |_{p^{-1}(x)} : p^{-1}(x) \ra ( H_d^r , H_d^s )  \ , \ 
t \in {(\rhors)_\eps} 
\right\} \ .
\]
\noindent If $t \in {(\rhors)_\eps}$, then by (\ref{rel_inv}) we find that the norm function $n_t \in C(\Omega)$: $n_t (\omega) := \left\| t (\omega) \right\|$ is constant on $p^{-1}(x)$, $x \in X$ (in fact, $n_t (\omega) = n_t (\omega \cdot g)$, $g \in G$); so that, $n_t \in C(X)$ and $\left\| t \right\| = \left\| n_t \right\|$. In other terms, the vector field $(t_x)_{x \in X} \in \prod_x \rho^{r,s}_x$ has a continuous norm function. We define the set $\Theta_{r,s} := \left\{ (t_x) \in \prod_x \rho^{r,s}_x \ : \ t \in {(\rhors)_\eps} \right\}$; by construction, $\Theta_{r,s}$ defines a unique continuous field of Banach spaces (\cite[Prop.10.2.3]{Dix}); in such a way, we interpret ${(\rhors)_\eps}$ as the module of continuous sections of the corresponding bundle of Banach spaces.

Let us now prove that $\rho^{r,s}_x \simeq ( H_d^r , H_d^s )_{L_x}$ for some closed Lie group $L_x \subseteq G_0$. We consider a set of the type $\left\{ \omega_x \in p^{-1} (x) \right\}_{x \in X} \subseteq \Omega$, and define $L_x := \left\{ u \in G_0 \ : \ \omega_x \circ \alpha_{ x , u } = \omega_x \right\}$. If $t \in {(\rhors)_\eps}$, $u \in L_x$, then (\ref{rel_inv}) implies $t (\omega_x) = t (\omega_x \circ \alpha_{ x , u } ) = u_s \cdot t(\omega_x) \cdot u_r^*$. Thus, $t (\omega_x) \in ( H_d^r , H_d^s )_{L_x}$. We define the linear map 
\[
\pi_x : \rho^{r,s}_x \ra ( H_d^r , H_d^s )_{L_x} \ \ , \ \ 
\pi_x (t_x) := t (\omega_x) \ .
\]
\noindent By (\ref{rel_inv}), we find that $\pi_x$ is isometric. Furthermore, let $t_0 \in ( H_d^r , H_d^s )_{L_x}$; then, a vector field $t_1 \in (p_* \mE^r , p_* \mE^s) |_{p^{-1}(x)} \simeq C ( p^{-1} (x)) \otimes (H_d^r , H_d^s)$ is defined by $t_1 (\omega_x \cdot u) := \wa u (t_0) = u_s t_0 u_r^*$, $u \in G_0$. By construction, $t_1$ is $G_0$-invariant w.r.t. the action (\ref{def_gact}). Let $W \subseteq X$ be a closed set trivialiazing $\mG$ as in Lemma \ref{lem_gcut}; by the Tietze theorem \cite[Lemma 1.6.3]{Ati}, there exists a vector field $t_2 \in (p_* \mE^r_W , p_* \mE^s_W)_{G_0}$ extending $t_1$. By (\ref{eq_def_ra}), there exists a vector field $t_3 \in (p_* \mE^r , p_* \mE^s)_G$ extending $t_2$, so that $\pi_x (t_3) = t_0$. Thus $\pi_x$ is also surjective, and the first assertion is proved. The second assertion trivially follows from the first one.
\end{proof}

\section{Special endomorphisms induce dual actions.}
\label{cp_spec}

\subsection{The main result.}

The following basic result is a direct consequence of Lemma \ref{lem_da}, Lemma \ref{dim_perm}, and Thm.\ref{soro_field}. It supplies a bundle descprition of the category $\wa \rho_\eps$ associated with a weakly special endomorphism. The hypothesis $c_1(\rho) \in \bN$ is made just in order to simplify the exposition (see the remarks following Def.\ref{dim_pm}).

\begin{thm}
\label{sue_action}
Let $\rho$ be a weakly special endomorphism of a unital \sC algebra $\mA$, with class $c(\rho) = d \oplus c_1(\rho)$, $d \in \bN$. Then, for every rank $d$ vector bundle $\mE \ra \spzro$ with first Chern class $c_1(\rho)$, there is a dual action $$\mu : \wa \mSUE \ra \wa \rho \subset \tend \ .$$ Moreover, $\soro$ is a continuous bundle of \sC algebras over $\spzro$, with fibres isomorphic to $\mO_{G^x}$, $x \in \spzro$, where $G^x \subseteq \sud$ is a compact Lie group unique up to conjugacy in $\sud$. Finally, there is an inclusion $\mO_{\mSUE} \hra \soro$ of \sC algebra bundles.
\end{thm}

\begin{cor}
For every $r,s \in \bN$, the symmetry intertwiners spaces $(\rhors)_\eps$ are continuous fields of Banach spaces over $\spzro$, with fibres $(H_d^r , H_d^s)_{G^x}$, $x \in X$. Moreover, with the notation of Prop.\ref{rhoine}, there is an isomorphism $\wa \rho_\eps \simeq p_* \wa \mSUE$, i.e. $(\rhors)_\eps \simeq ( p_* \mE^r , p_* \mE^s )_{\mSUE}$, $r,s \in \bN$.
\end{cor}

\begin{ex} {\em 
\label{ex_bl}
In the case of canonical endomorphisms considered in \cite{BL04}, it turns out that every $(\rhors)_\eps$ is a trivial field for a fixed compact group $G$: $(\rhors)_\eps \simeq \zro \otimes (H_d^r , H_d^s)_G$. The fibre $(H_d^r , H_d^s)_G$ corresponds to the intertwiners space $(\rho^r , \rho^s)_\bC$ (about the previous notation, see Thm.4.10 of the above-cited reference). }
\end{ex}

\subsection{Final Remarks.}


Let $\rho \in \tend$ be a quasi-special endomorphism; by using Thm.\ref{sue_action}, it is possible to construct the crossed product $\mA \rtimes_\mu \wa {\mSUE}$. By Lemma \ref{lem38}, the fixed-point algebra w.r.t. the action $$\alpha : \mSUE \ra {\bf aut}_{\spzro} (\mA \rtimes_\mu \wa {\mSUE})$$ is $\mA$. In order to establish an isomorphism of tensor \sC categories $\wa \rho_\eps \simeq \wa G$ for some closed group $G \subseteq \mSUE$, we should construct a crossed product $\mc$ with the condition $\mA' \cap (\mc) = \mA \cap \mA'$, as by Prop.\ref{cor_dual_1}. Such a crossed product should be obtained as a suitable quotient of $\mA \rtimes_\mu \wa {\mSUE}$.

At this purpose, let us  recall the original construction by Doplicher and Roberts for the trivial-centre case ($\mA \cap \mA' = \bC1$ $\Rightarrow$ $\mE = \bC^d$). As a first step, the crossed product $\mB := \mA \rtimes_\mu \wa{\sud}$ is considered; then, the following quotient is defined:
\[
\mF := \mB \ / \ (C_\omega(\Omega) \ \mB ) \ \ ,
\]
\noindent where $\mC := C(\Omega)$ is the centre of $\mB$, $\omega \in \Omega$, $C_\omega(\Omega) := \left\{ c \in \mC : c(\omega) = 0 \right\}$; we denote by $\eta : \mB \ra \mF$ the associated epimorphism. Now, by construction there is an action $\alpha : \sud \ra {\bf aut} (\mB)$; we define
\[
G := \left\{ g \in \sud : \omega \circ \alpha_g = \omega \right\} \ .
\]
\noindent It is clear that $C_\omega (\Omega) \mB$ is $G$-stable, thus we may introduce the action $\beta : G \ra {\bf aut} \mF$, $\beta_g \circ \eta = \eta \circ \alpha_g$, $g \in G$.  At this point, an isomorphism $\mF \simeq \mc$ is established, and it turns out that $\eta (\mA)' \cap \mF = \bC 1$ (see \cite[Thm.4.1]{DR89A}).

In order to perform the analogous construction for $\mA \cap \mA' \neq \bC 1$, the pure state $\omega$ has to be replaced with a $C(X)$-epimorphism $\varphi : \mC \ra C(X)$. In such a way, we may define $\mF := \mB  /  ( \ker \varphi \ \mB )$, and conclude $\mF \simeq \mc$ in the same way as Doplicher and Roberts. In the case considered by Baumg\"artel and Lled\'o (\cite{BL04}), it is easy to construct such a $\varphi$: in fact, by using \cite[Lemma 5.1]{DR88} it can be proved that $\mC \simeq C(X) \otimes \mC'$, where $\mC'$ is the \sC algebra generated by elements of the type $t^* \psi$, $t \in (\iota , \rho^r)_\bC$, $\psi \in H_d^r$, $r \in \bN$ (see Ex.\ref{ex_bl}). So that, we may define $\varphi := id_X  \otimes \omega'$, where $\omega'$ belongs to the spectrum of $\mC'$.

In the general case, the existence of $\varphi$ is not trivial to prove. Moreover, the group $G$ reconstructed in this way depends on the choice of $\mE$ (see \cite[Cor.4.18, Ex.4.7]{Vas04} for examples in the case $G = \mSUE$). These questions, relative to existence and unicity problems, will be approached in a forthcoming paper.


\

{\bf Acknowledgments.} The author would like to thank S. Doplicher, for having proposed (and successively supported with improvements and encouragements) the initial idea of this work for his PhD Thesis.


\end{document}